\documentclass[reqno,12pt]{amsart}

\usepackage{amsmath,amsthm,amscd,amssymb}

\usepackage{latexsym}

\usepackage{mathrsfs}

\newcommand{\bbN}{{\mathbb{N}}}

\newcommand{\bbR}{{\mathbb{R}}}

\newcommand{\bbZ}{{\mathbb{Z}}}

\newcommand{\bbC}{{\mathbb{C}}}

\newcommand{\calF}{{\mathcal F}}

\newcommand{\calD}{{\mathcal D}}

\newcommand{\calU}{{\mathcal U}}

\newcommand{\calX}{{\mathcal X}}

\newcommand{\calW}{{\mathcal W}}

\newcommand{\calV}{{\mathcal V}}

\newcommand{\bfd} {{\mathbf d}}

\newcommand{\sone}{{\bf 1}}

\newcommand{\szer}{{\bf 0}}

\newcommand{\se}[1]{\{#1(n)\}}

\newcommand{\sse}[2]{\{#1(n)+#2(n)\}}

\newcommand{\an}{\{a(n)\}}

\newcommand{\bn}{\{b(n)\}}

\newcommand{\Jab}{J \left( \an,\bn \right)}

\newcommand{\mean}[1]{\left \langle #1 \right \rangle}

\newcommand{\no}{\nonumber}

\newcommand{\f}{\frac}

\newcommand{\ti}{\tilde  }

\newcommand{\ac}{\text{\rm{ac}}}

\newcommand{\singc}{\text{\rm{sc}}}

\newcommand{\pp}{\text{\rm{pp}}}

\newcommand{\beq}{\begin{equation}}

\newcommand{\eeq}{\end{equation}}

\newcommand{\ba}{\begin{align}}

\newcommand{\ea}{\end{align}}

\newcommand{\veps}{\varepsilon}

\newcommand{\vphi}{\varphi}

\newcommand{\alphac}{\alpha\textrm{c}}

\newcommand{\alphas}{\alpha\textrm{s}}

\allowdisplaybreaks

\numberwithin{equation}{section}

\newtheorem{theorem}{Theorem}[section]

\newtheorem{proposition}[theorem]{Proposition}

\newtheorem{lemma}[theorem]{Lemma}

\newtheorem{corollary}[theorem]{Corollary}

\theoremstyle{definition}

\theoremstyle{remark}

\newtheorem*{remark}{Remark}

\newcommand{\abs}[1]{\lvert#1\rvert}

\begin{document}

\title[Stability of Spectral Types]{Stability of Spectral Types for 
Jacobi Matrices Under Decaying Random Perturbations}


\author[J.~Breuer and Y.~Last]{Jonathan Breuer$^{1,2}$ and Yoram Last$^{1,3}$}

\thanks{$^1$ Institute of Mathematics, The Hebrew University,
91904 Jerusalem, Israel.}
\thanks{$^2$ E-mail: jbreuer@math.huji.ac.il.}
\thanks{$^3$ E-mail: ylast@math.huji.ac.il.}

\date{December 4, 2006}

\begin{abstract}
We study stability of spectral types for 
semi-infinite self-adjoint tridiagonal matrices under
random decaying perturbations. We 
show that absolutely continuous spectrum associated with bounded 
eigenfunctions is stable under Hilbert-Schmidt random perturbations. We 
also obtain some results for singular spectral types. 
\end{abstract}

\maketitle

\section{Introduction}
In this paper we study semi-infinite Jacobi matrices of the form
\beq \label{jacobi}
J \left( \{a(n)\}_{n=1}^\infty,\{b(n)\}_{n=1}^\infty \right)= \left( 
\begin{array}{ccccc}
b(1)    & a(1) & 0      & 0      & \ldots \\
a(1)    & b(2) & a(2)    & 0      & \ldots \\
0      & a(2) & b(3)    & a(3)    & \ddots \\
\vdots & \ddots   & \ddots & \ddots & \ddots \\
\end{array} \right)
\eeq
with \beq \no  b(n) \in \bbR, \ a(n)>0, \eeq
as operators
on $\ell^2(\bbZ_+=\{1,2, \ldots \})$. We shall assume throughout that 
$\Jab$ is self-adjoint. For this to be true, $\Sigma_{n=1}^\infty a(n)^{-1}=\infty$ 
suffices \cite{berez}. In fact, we need a somewhat stronger restriction on the growth of 
$\se{a}$ (see \eqref{growth-restrict} below). 
 
Such operators are a natural generalization of discrete Schr\"odinger 
operators on the half line. In particular, the discrete Laplacian on $\ell^2(\bbZ_+)$
can be described with the help of the constant 
sequences $\sone\equiv\se{a^\circ}$, $\szer\equiv\se{b^\circ}$, where
$a^\circ(n) = 1$ and $b^\circ(n)= 0$ for all $n\in\bbZ_+$, so that 
\beq \no
\Delta=J \left( \sone, \szer \right).
\eeq

 From the fact that the vector 
\beq \no 
\delta_1 \equiv \left( \begin{array}{c}
1 \\
0 \\
0 \\
\vdots \\
\end{array} \right)
\eeq 
is a cyclic vector for 
$\Jab$, it follows (\cite{reed-simon1})
that there exists a measure $\mu$, which coincides with the spectral
measure of the vector $\delta_1$, so that $\Jab$ is unitarily
equivalent to the operator of multiplication by the parameter on 
$L^2(\bbR,d\mu)$. $\mu$ decomposes as 
\begin{displaymath}
\mu=\mu_{\ac}+\mu_{\singc}+\mu_{\pp},
\end{displaymath} 
where $\mu_{\ac}$ is the part of $\mu$
that is absolutely continuous with respect to the Lebesgue measure,
$\mu_{\singc}$ is a continuous measure that is singular with 
respect to the Lebesgue measure, and $\mu_{\pp}$ is a pure point
measure.

We want to investigate the stability of certain continuity properties
of $\mu$ under a decaying random perturbation of
$\Jab$. The first part of the paper deals with the stability of 
the essential support of the absolutely continuous spectrum. In the 
second part, we restrict the discussion to the case
$\se{a}=\sone$ (the discrete Schr\"odinger case) and deal
with the more delicate singular spectral types. In 
both cases, a principal tool in the analysis is the connection 
between properties of the spectral measure and the behavior at infinity 
of solutions of the difference equation
\beq \label{difference} 
a(n)\vphi(n+1)+a(n-1)\vphi(n-1)+b(n)\vphi(n)=E\vphi(n) 
\eeq 
for fixed $E\in\bbR$ and $n\geq 1$ (we set $a(0)=1$). Such a difference 
equation can be 
regarded 
as an initial value problem, which makes it natural to introduce the 
\emph{single-step transfer matrices}:
\beq \label{transfer}
S^E(n)=\left( \begin{array}{cc}
\frac{E-b(n)}{a(n)} & -\frac{a(n-1)}{a(n)} \\
1 & 0 \end{array} \right), \quad n\geq1 ,
\eeq 
that satisfy 
\beq \no \left( 
\begin{array}{c} \vphi(n+1) \\
\vphi(n) \end{array} \right)=S^E(n)\left(
\begin{array}{c} \vphi(n) \\
\vphi(n-1) \end{array} \right) 
\eeq
for any $\{\vphi(n)\}_{n=0}^\infty$ that solves \eqref{difference}. Thus, if 
we 
denote 
\beq \no 
\vec{\vphi}(n)=\left( \begin{array}{c} \vphi(n+1) \\
\vphi(n) \end{array} \right)
\eeq
and $T^E(n)\equiv S^E(n)\cdot \ldots \cdot S^E(1)$,
then 
\beq \label{n-transfer} 
\vec{\vphi}(n)=T^E(n) \vec{\vphi}(0).
\eeq
\sloppy

The \emph{essential support} of an absolutely continuous measure $\nu$ 
on $\bbR$ is the equivalence class $\Sigma_\ac(\nu)$ of sets $A 
\subseteq \bbR$ such that $\nu$ is supported on $A$ and that the 
restriction of Lebesgue measure to $A$ is absolutely continuous w.r.t.\ 
$\nu$.
We shall use $\Sigma_\ac(\se{a},\se{b})$ to denote the essential 
support of $\mu_\ac$ and refer to it as the essential support of the 
absolutely continuous spectrum of $\Jab$.

Over the past decade, there has been a significant amount of work done (see, e.g., 
\cite{christ-kiselev, christ-kiselev-b, christ-kiselev-remling, deift-killip,
killip, kiselev1, kiselev2, remling}), in 
the area of one-dimensional Schr\"odinger operators,
towards determining conditions on a perturbing potential $\se{\ti{b}}$ 
ensuring that
\beq \label{stability} 
\Sigma_\ac(\sone,\se{b})=\Sigma_\ac(\sone,\sse{b}{\ti{b}}). 
\eeq 
That such an equality exists for any $\se{\ti{b}}\in\ell^1$ is a well 
known 
result
 from scattering theory \cite[Chapter~XI.3]{reed-simon3}. 
For general $\se{b}$, this is the best there is at present, in terms of 
sheer $\ell^p$ properties of the perturbation.  For $\se{b}=\szer$, 
however,
it has been proven by Deift-Killip \cite{deift-killip}
that \eqref{stability} holds for $\se{\ti{b}}$ merely in $\ell^2$. This 
result has been later extended by Killip \cite{killip} to include
any periodic $\se{b}$.
For arbitrary background potentials $\se{b}$, it 
has been conjectured 
by Kiselev-Last-Simon \cite{kls} that an $\ell^2$ 
perturbation does not change the essential support of the absolutely 
continuous spectrum.
For a perturbation of the off-diagonal entries as well as 
the diagonal entries, Killip-Simon \cite{killip-simon} have shown that if 
$\se{\ti{a}},\se{\ti{b}}\in\ell^2$, then   
\beq \label{stability1}
\Sigma_\ac(\sone,\szer)=\Sigma_\ac(\sone+\se{\ti{a}},\se{\ti{b}}).
\eeq

Our first result deals with the preservation of $\Sigma_\ac(\se{a},\se{b})$ for 
general $\se{b}$ and $\se{a}$ obeying 
\beq \label{growth-restrict}
\limsup_{L\rightarrow \infty} \frac{1}{L}\sum_{n=1}^L a(n)^{-1}>0
\eeq 
under a \emph{random} decaying 
perturbation of both the diagonal and off-diagonal entries.  
For a measurable set $B\subseteq\bbR$, $\Sigma_\ac\cap B$ denotes
the equivalence class of sets $A\cap B$ such that $A\in\Sigma_\ac$.

\begin{theorem} \label{ac_stability} 
Let $\Jab$ be a Jacobi 
matrix such that $\se{a}$ obeys \eqref{growth-restrict}, and let $\ti{a}_\omega(n):\Omega\rightarrow\bbR$ and 
$\ti{b}_\omega(n):\Omega\rightarrow\bbR$ $(n\geq1)$ be 
two sequences of independent random variables with zero mean, defined 
over a 
probability space $(\Omega,\mathcal{F},P)$. Assume that there 
exists a $\delta>0$, for which 
\beq \label{condition_on_alpha}
\delta^{-1}>\frac{a(n)}{a(n)+\ti{a}_\omega(n)}>\delta
\eeq
for every $n$ and $\omega\in\Omega$. Let $J_0=\Jab$ and 
$$
J_\omega=J(\sse{a}{\ti{a}_\omega},\sse{b}{\ti{b}_\omega}) .
$$
Then, for a.e.\ $\omega$,
\beq \label{ac_conclusion} 
\Sigma_\ac(J_0)\cap\Gamma=\Sigma_\ac(J_\omega)\cap\Gamma ,
\eeq
where $\Gamma$ is the set of all $E\in \bbR$ for which
\beq \label{Gamma-def}
\sum_{n=1}^\infty
\left( \mean{\ti{a}_\omega(n)^4}^{1/2}
+ \mean{\ti{b}_\omega(n)^2} \right)
\left((a(n)+1)t^E(n)\right)^4 <\infty ,
\eeq
where we denote $\mean{f_\omega} \equiv \int_\Omega f_\omega dP(\omega)$
for any measurable function $f_\omega$ of $\omega$
and \mbox{$t^E(n) \equiv \parallel T^E(n)\parallel$}
is the norm of the $n$'th transfer matrix corresponding to $J_0$.
\end{theorem}

We note that Kaluzhny-Last \cite{kala} recently studied Jacobi matrices
of the form $J(\sse{a}{\ti{a}_\omega},\sse{b}{\ti{b}_\omega})$, where
$\se{a}-\sone$ and $\se{b}$ are decaying sequences of bounded variation
and $\se{\ti{a}_\omega}$, $\se{\ti{b}_\omega}$ are as in Theorem \ref{ac_stability}
and obey
$$
\sum_{n=1}^\infty
\left( \mean{\ti{a}_\omega(n)^2}
+ \mean{\ti{b}_\omega(n)^2} \right)
<\infty .
$$
They show that, with probability one, such operators have purely
absolutely continuous spectrum on $(-2,2)$ and moreover, this purity
of the absolutely continuous spectrum is stable under changing any finite
number of entries in the Jacobi matrices. Since the unperturbed
$J(\se{a},\se{b})$ is known (see, e.g., \cite{bounded-eig}) in this
case to have purely absolutely continuous spectrum on $(-2,2)$
with $\{t^E(n)\}_{n=1}^\infty$ being a bounded sequence for every
$E\in (-2,2)$, we see that a part of their result, namely, the
fact that $\Sigma_\ac(J_0)=\Sigma_\ac(J_\omega)$, can be recovered
as a special case of Theorem \ref{ac_stability}.

To further elucidate Theorem \ref{ac_stability},
consider the case $a(n)=1$, $\ti{a}_\omega(n)=0$. The 
condition defining $\Gamma$ translates into an $\ell^2$ type condition
on the perturbation when one studies energies for which the transfer 
matrices
are bounded: For a given background potential $\se{b}$, denote
\beq \no 
\Gamma_0\equiv\Gamma_0(\se{b})=\{E\in\bbR \mid t^E(n) \textrm{ is bounded}\}. 
\eeq 
Then it follows from the theory of subordinacy \cite{subord} (also see
\cite{bounded-eig})
that there exists a set $A\in\Sigma_\ac(\sone, \se{b})$ for which 
$\Gamma_0\subseteq A$. From Theorem~\ref{ac_stability}, it 
follows that
\begin{corollary} \label{corollary-ac}
Assume that 
$$
\sum_{n=1}^\infty \mean{\ti{b}_\omega(n)^2}<\infty .
$$
Then, for a.e.\ $\omega$,
\beq \no
\Sigma_\ac ( \sone,\se{b} )\cap\Gamma_0= 
\Sigma_\ac ( \sone,\sse{b}{\ti{b}_\omega} )\cap\Gamma_0.
\eeq
\end{corollary}
Corollary~\ref{corollary-ac} constitutes some progress towards a random version 
of the above mentioned conjecture of Kiselev-Last-Simon \cite{kls}.
Whether actually $\Gamma_0(\se{b})\in\Sigma_\ac(\se{b})$ for any $\se{b}$
is a long standing open problem.
For some related work, see Maslov-Molchanov-Gordon \cite{schrodinger-conjecture}.

The question of stability of \emph{singular} spectral types has received 
much less attention than the one concerning $\Sigma_\ac$. One of the 
reasons for this is the fact that singular spectral types are not stable 
even under rank one perturbations (see \cite{djls}). 
One may, however, bypass this problem by using an idea 
of Del-Rio-Simon-Stolz \cite{delrio} to consider the union of spectral 
supports over the different boundary conditions. This provides a unified 
approach for the different spectral types, in 
that spectral stability is obtained for any compactly supported 
perturbation (see \cite{delrio}). Kiselev-Last-Simon \cite{kls} have 
modified 
and extended this approach, via 
the theory of subordinacy, to deal with the classification of spectral types 
according to the singularity/continuity of the spectral measure 
w.r.t.\ $\alpha$-dimensional Hausdorff measures. 
In our definitions, we follow their general methodology. 

While it is possible, using the methods developed below, to deal with the 
general Jacobi case, we restrict the discussion to the case of diagonal
perturbations of discrete Schr\"odinger operators. We take this approach 
in order to avoid technical difficulties which may obscure the main 
argument.  
Thus, for fixed $E\in\bbR$, we shall be looking at properties of solutions 
of the equations 
\beq \label{difference-schrodinger}
\vphi(n+1)+\vphi(n-1)+b(n)\vphi(n)=E\vphi(n)
\eeq
for $n\geq2$,
\beq \label{boundary-potential}
\vphi(2)+(b(1)-\tan(\theta))\vphi(1)=E\vphi(1)
\eeq
for $-\frac{\pi}{2}<\theta<\frac{\pi}{2}$.
Such sequences are obviously eigenvectors (not necessarily in $\ell^2$) of 
the infinite matrix
\beq \label{Htheta}
H_{\theta}=\left( \begin{array}{ccccc}
b(1)-\tan(\theta)    & 1 & 0      & 0      & \ldots \\
1    & b(2) & 1    & 0      & \ldots \\
0      & 1 & b(3)    & 1    & \ddots \\
\vdots & \ddots   & \ddots & \ddots & \ddots \\
\end{array} \right).
\eeq

We denote by $\vphi^E_{1,\theta}(n)$ the solution to 
\eqref{difference-schrodinger}, \eqref{boundary-potential}, normalized by 
\beq \label{phi1theta} \vphi^E_{1,\theta}(1)=\cos(\theta).
\eeq
We also include the case $\theta=-\pi/2$, for which
\eqref{boundary-potential} and \eqref{phi1theta}
are replaced by $\vphi^E_{1, -\pi/2}(1)=0$,
$\vphi^E_{1, -\pi/2}(2)=1$.
We shall use the notation
\beq \label{phi2theta} 
\vphi^E_{2,\theta}\equiv \vphi^E_{1, \theta-\pi/2}.
\eeq
\begin{remark}
One may define $\vphi^E_{1,\theta}$ by referring only to 
\eqref{difference-schrodinger} (for $n\geq1$) and using 
$\vphi^E_{1,\theta}(0)=-\sin(\theta)$, $\vphi^E_{1,\theta}(1)=\cos(\theta)$.
This way $\vphi^E_{1,\theta}$ is more naturally defined on
$[-\frac{\pi}{2},\frac{\pi}{2})$, without anything special
for $\theta=-\pi/2$.
\end{remark}
A basic object in the theory of subordinacy is the $L$'th norm (for 
$L>0$) 
of a function $f:\bbZ_+\rightarrow\bbC$,
\beq \parallel f 
\parallel_L^2  \equiv \sum _{n=1}^{\lfloor L 
\rfloor} \abs{f(n)}^2+(L-\lfloor L \rfloor)\abs{f(\lfloor L \rfloor + 1)}^2 , 
\eeq
where $\lfloor \,\cdot\, \rfloor$ denotes integer part.
For a given 
$E\in\bbR$, $\theta\in[-\frac{\pi}{2},\frac{\pi}{2})$, 
$\vphi^E_{1,\theta}$ is called \emph{subordinate} if 
\beq \label{subordinacy} 
\lim_{L\rightarrow\infty}\f{\parallel 
\vphi^E_{1,\theta} \parallel_L}{\parallel \vphi^E_{2,\theta} \parallel_L}=0
.
\eeq
It is clear that a subordinate solution does not necessarily exist for 
every $E$, but whenever it does, it is unique. We denote the $\theta$ 
for which $\vphi^E_{1,\theta}$ is subordinate, if it exists,
by $\theta(E)$. One may 
decompose $\bbR$ into three disjoint sets:
\beq \no \Sigma_{\pp}\equiv\{E\in\bbR\mid 
\theta(E)\text{ exists and }\vphi^E_{1,\theta(E)}\in\ell^2\}
\eeq
\beq \no \Sigma_{\singc}\equiv\{E\in\bbR\mid 
\theta(E)\text{ exists and }\vphi^E_{1,\theta(E)}\not\in\ell^2\}
\eeq
\beq \no \bbR\setminus (\Sigma_{\pp} \cup \Sigma_{\singc})
\eeq
What makes the discussion of stability of singular spectral types 
interesting is the fact (see, e.g., \cite{kls}) that these three sets have 
the following spectral interpretation:
\begin{itemize}
\item $\Sigma_{\pp}=\cup_{\theta}\ti{\sigma}_{\pp}(H_{\theta})$, 
where $\ti{\sigma}_{\pp}(H_{\theta})$ is the set of eigenvalues of
$H_{\theta}$. 

\item For any $\theta$, 
$\mu_{\theta,\singc}(\cdot)=\mu_{\theta}(\Sigma_{\singc}\cap\cdot)$ 
and any other set $A$ with this property equals $\Sigma_{\singc}$ up to 
a set of Lebesgue measure zero. 

\item $\Sigma_{\ac}\ni\bbR 
\setminus (\Sigma_{\pp} \cup \Sigma_{\singc}).$ 
\end{itemize}
The above sets are clearly independent of $\theta$ and stable under 
compactly supported perturbations.

The Jitomirskaya-Last extension of subordinacy theory \cite{jit-last} 
makes it possible to 
investigate the stability of Hausdorff-dimensional properties of the 
spectral measure. It follows from their analysis that for any 
$\alpha\in(0,1]$, 
there exist sets $\Sigma_{\alphac}\subseteq\bbR$ and 
$\Sigma_{\alphas}\subseteq\bbR$ such that for any $\theta$, 
\beq  \label{alpha_supports}
\mu_{\theta,\alphac}=\mu_\theta(\Sigma_{\alphac}\cap\cdot), \quad
\mu_{\theta,\alphas}=\mu_\theta(\Sigma_{\alphas}\cap\cdot)
\eeq
where $\mu_{\theta,\alphac}$ is the part of $\mu_\theta$ that is 
continuous with respect to the $\alpha$-dimensional Hausdorff measure, and 
$\mu_{\theta,\alphas}$ is the part which is singular with respect to it. 
(For the study of decompositions of a measure w.r.t.\ 
dimensional Hausdorff measures and for the significance of 
this analysis to quantum mechanics, see, for example, 
\cite{hausdorff} and references therein.) For any 
$\alpha\in(0,1]$,
$\Sigma_{\alphas}\subseteq\Sigma_{\singc}\cup\Sigma_{\pp}$, and for any $E\in
\Sigma_{\singc}$, whether $E\in 
\Sigma_{\alphas}$ or not, depends on the decay of the subordinate solution 
at infinity: 
\beq \no
E\in \Sigma_{\alphas}
\eeq
if and only if
\beq \label{singularity_condition}
\liminf_{L \rightarrow \infty} \frac{\parallel \vphi^E_{1,\theta 
(E)}\parallel _{L}}{\parallel \vphi^E_{2,\theta (E)}\parallel _{L}^
{\ti{\beta}(\alpha) }}=0
\eeq
where  $\ti{\beta}(\alpha) =\frac{\alpha }{2-\alpha}$ (see \cite{jit-last}).

The discussion above motivates the following definition of 
\cite{kls}:
Let $E\in\Sigma_{\singc}$. Define 
\beq
\beta(E)=\liminf_{L \rightarrow \infty} 
\frac{\ln \parallel \vphi^E_{1,\theta(E)} \parallel _L}{\ln \parallel \vphi^E_{2,\theta(E)} \parallel _L}.
\eeq
For any $E$ with $\beta(E)>0$, we also define 
\beq \label{eta(E)}
\eta(E)=\frac{1-\beta(E)}{\beta(E)}.
\eeq
Again, it is clear that the sets $\Sigma_{\alphas}$ and $\Sigma_{\alphac}$ 
and the parameter $\beta(E)$ (where it is defined) are stable under 
compactly supported perturbations.
To obtain more, one needs a regularity condition on the energy:
Following Kiselev-Last-Simon \cite{kls}, we shall call
an energy $E$ \emph{regular} if for some $\theta$ and all 
$\veps>0$, we have 
\beq \no
\parallel \vphi^E_{1,\theta} \parallel _L<C_{\veps}L^{\frac{1}{2}+\veps}.
\eeq
Since almost every energy is regular both with respect to each
$\mu_\theta$ (see \cite{berez}) and 
(by spectral averaging---see Theorem~1.8 in \cite{rankone}) with respect 
to Lebesgue measure, the demand that 
energies be regular is not a severe restriction.

Let 
\beq \label{lambda_naught}
\Lambda_0=\{E\in\Sigma_{\singc}\mid E \textrm{ is regular and } 
\beta(E)>0\} .
\eeq
For deterministic perturbations, Kiselev-Last-Simon \cite{kls} have 
shown 
that,
for any $0<\alpha<1$,
\beq \no
\ti{\Lambda}\cap\Sigma_{\alphac} \left( \sone,\se{b} \right)
\subseteq\Sigma_{\alphac} \left( \sone,\sse{b}{\ti{b}} \right),
\eeq
\beq \no 
\ti{\Lambda}\cap\Sigma_{\alphas}\left( \sone,\se{b} \right) 
\subseteq\Sigma_{\alphas} \left( \sone,\sse{b}{\ti{b}} \right),
\eeq
where $\se{b}$ is any background potential, and 
\beq \label{Lambda_tilde}
\ti{\Lambda}=\{E\in\Lambda_0\mid  \vert \ti{b}(n) \vert<Cn^{-\gamma}
\textrm{ for some } \gamma>\eta(E)+1\}.
\eeq

For random potentials we show
\begin{theorem} \label{singular_stability}
Let $\se{\ti{b}_\omega}$ be a sequence of independent
real-valued random
variables with zero mean on $\Omega$. For any $E\in\Lambda_0$, 
$\ti{\eta}>0$ 
and $n\geq1$, let
\beq \label{definition-rnE}
r_{\ti{\eta}}^E(n)\equiv|\vphi^E_{1,\theta(E)}(n)|^4n^{2\ti{\eta}}+|\vphi^E_{2,\theta(E)}(n)|^4,
\eeq
and let
\beq \label{Lambda}
\Lambda=\left\{E\in\Lambda_0\;\Bigg|\;
\sum_{n=1}^\infty \left( r_{\ti{\eta}}^E(n) \mean{\ti{b}_\omega(n)^2} \right)<\infty \textrm{ for some } 
\ti{\eta}>\eta(E)\right\}.
\eeq
Then, for any $0<\alpha<1$ and any fixed measure $\nu$ on $\bbR$,
for a.e.\ $\omega$,
\beq \no
\Lambda\cap\Sigma_{\alphac} ( \sone, \se{b} )
\subseteq \Sigma_{\alphac} ( \sone, \sse{b}{\ti{b}_\omega} ),
\eeq
\beq
\Lambda\cap\Sigma_{\alphas} ( \sone, \se{b} )
\subseteq \Sigma_{\alphas} ( \sone, \sse{b}{\ti{b}_\omega} ),
\eeq
where the inclusion is up to a set of $\nu$-measure zero. 
\end{theorem}
Theorems \ref{ac_stability} and 
\ref{singular_stability} have the common feature of the appearance of the 
4th power of the norms of the transfer matrices (in \ref{singular_stability}, 
see the definition of $r_{\ti{\eta}}^E(n)$). The reason for this is 
that our basic tool is a random variation of parameters, where the 
perturbing potential is 
coupled to the square of the transfer matrices and thus, when estimating 
the variance of the perturbation, the 4th power enters the picture.
For examples where the pointwise behavior of the solutions to 
\eqref{difference} is known, this does not constitute a problem. 
One such example is the class of bounded sparse potentials 
studied by Zlato\v s \cite{zlatos}. For this class, one has stability of 
$\Sigma_{\alphac}$ and $\Sigma_{\alphas}$ under random perturbations 
decaying like $n^{-\gamma}$ for $\gamma>\eta(E)+\frac{1}{2}$ (compare 
with 
$\gamma>\eta(E)+1$ in \eqref{Lambda_tilde}). 

In light of these remarks, the general question of the pointwise 
behavior of the solutions of \eqref{difference} is one that arises naturally 
in connection with the results presented here. The more famous question of 
whether or not there is almost-everywhere boundedness of solutions
with respect to the 
absolutely continuous part of the spectral measure is only one facet of 
this general problem.

The rest of this paper is organized as follows. Section 2 covers some 
preliminaries---especially a useful characterization of $\Sigma_{\ac}$ due 
to Last-Simon \cite{last-simon} and a variation on a classic 
theorem concerning the almost 
everywhere convergence of random series with convergent variances. 
In Section 3 we introduce the main idea behind our analysis.
We formulate and prove two different (but similar) lemmas
which are central to the proofs of our two main theorems. 
These theorems are proved in Section 4. 
Section 5 has the worked out 
application of Theorem~\ref{singular_stability} to the
above mentioned sparse potentials of
Zlato\v s \cite{zlatos}.

This research was supported in part
by The Israel Science Foundation (Grant No.\ 188/02)
and by Grant No.\ 2002068 from the United States-Israel
Binational Science Foundation (BSF), Jerusalem, Israel.


\section{Preliminaries}

As explained in the introduction, we want to exploit the connection 
between spectral properties of the operator $\Jab$ and 
the asymptotic properties of the solutions to the corresponding difference 
equation. That is, we want to compare the asymptotic properties of the 
solutions to the difference equation corresponding to the basic 
operator, with those of the solutions to the equation corresponding to 
the perturbed one. In the singular continuous case we will `equate' the 
behavior at infinity of the perturbed and unperturbed solutions (in a 
sense to be precisely defined in Section 4). For the absolutely continuous
case, 
however, we need a little less. We rely on the following 
characterization of $\Sigma_{\ac}$ due to Last-Simon \cite{last-simon}:
\begin{proposition}[Last-Simon \cite{last-simon}] \label{last-simon} 
Let $\Jab$ be a self-adjoint Jacobi matrix such that $\se{a}$
satisfies \eqref{growth-restrict}, and let $T^E(n)$ be the 
corresponding transfer matrices defined by \eqref{transfer}-\eqref{n-transfer}.
Let $\{N_j\}_{j=1}^\infty$ be a sequence for which
\beq \no
\lim_{j \to \infty}  \frac{1}{N_j}\sum_{n=1}^{N_j}\frac{1}{a(n)} >0
\eeq
and let $\Sigma_{\ac} \equiv \Sigma_{\ac} \left( \se{a},\se{b} \right )$.
Then 
\beq \no
\left \{E\in\bbR\;\Bigg|\; \liminf_{j\to\infty} 
\frac{1}{N_j}\sum_{n=1}^{N_j} \parallel T^E(n) \parallel ^2<\infty \right \} \in \Sigma_{\ac}.
\eeq
\end{proposition}
\begin{remark}
This is actually a slight generalization of Theorem~1.1
of \cite{last-simon} to the general Jacobi 
case. Its proof is essentially the same as their proof.
\end{remark}  

The following are variants of a martingale inequality and 
convergence theorem which play a crucial role in 
the proofs of Lemmas~\ref{random_variation} and \ref{random_variation2}.

\begin{lemma} \label{kolmo-inequality}
Let $(\Omega, \calF, P)$ be a probability space and let 
$\se{x_\omega}$ be a sequence of independent random 
variables such that 
\beq \no 
\int_\Omega x_\omega(n)dP(\omega)\equiv\mean{x_\omega(n)}=0 
\eeq 
for all $n$.
Let 
\beq \no
z_\omega(n)=x_\omega(n)f_n(x_\omega(n+1),x_\omega(n+2),\ldots)
\eeq
where the $f_n$ are real-valued, measurable functions on $\bbR^\infty$.

Then, for any $N_1<N_2$ and $r\geq 0$,
\beq \label{k-inequality}
P\left( \left \{ \omega \;\bigg|\; \max_{N_1\leq n \leq N_2} |z_\omega(n)+\ldots+z_\omega(N_2)|> r \right\} \right)\leq 
\frac{\sum_{j=N_1}^{N_2}\mean{(z_\omega(j))^2}}{r^2}.
\eeq
\end{lemma}
\begin{proof}
Obviously, we may assume that $\mean{\left( z_\omega(n) \right)^2}<\infty$ for all $n$, since otherwise there is nothing to prove.
Denote
\beq \no 
Y_\omega(n)=\sum_{j=N_1}^{n-1}z_\omega(j) , \qquad
Q_\omega(n)=\sum_{j=n}^{N_2}z_\omega(j),
\eeq
and let
\beq \no
A_j=\left \{\omega \in \Omega \mid |Q_\omega(j)|>r;\; |Q_\omega(j+1)|,\ldots,|Q_\omega(N_2)|\leq 
r \right \}.
\eeq
Then, if $i<j$,
\beq \no
\mean{z_\omega(i) Q_\omega(j) \chi_j}= \mean{x_\omega(i)} \mean{f_i\left( x_\omega(i+1),\ldots \right) Q_\omega(j) \chi_j}=0
\eeq
where
\beq \no
\chi_j=\chi_{A_j}=\textrm{the characteristic function of }A_j,
\eeq
and thus,
\beq \no
\mean{\chi_j Y_\omega(j) Q_\omega(j)}=0
\eeq
so that 
\beq \no
\mean{\chi_j Q_\omega(j)^2} \leq \mean{\chi_j \left( Y_\omega(j)+Q_\omega(j) \right)^2}.
\eeq
Therefore
\beq \no
r^2 \mean{\chi_j} \leq \mean{\chi_j Q_\omega(j)^2} \leq \mean{\chi_j \left( Y_\omega(j)+Q_\omega(j) \right)^2}
\eeq
and
\beq \no
\begin{split}
r^2\sum_{j=N_1}^{N_2} \mean{\chi_j} &\leq 
\sum_{j=N_1}^{N_2} \mean{\chi_j Q_\omega(j)^2} 
\leq \sum_{j=N_1}^{N_2} \mean{\chi_j  \left( Y_\omega(j)+Q_\omega(j) \right)^2} \\
&=\sum_{j=N_1}^{N_2} \mean{\chi_j
\left( \sum_{l=N_1}^{N_2}z_\omega(l) \right)^2} \leq 
 \mean{ \left( \sum_{j=N_1}^{N_2}z_\omega(j) \right)^2} \\
&=\mean{\sum_{j=N_1}^{N_2}z_\omega(j)^2}
\end{split}
\eeq
where in the last equality we use
\beq \no
\mean{z_\omega(i)z_\omega(j)}=0 \quad \textrm{for }i\not=j.
\eeq
This ends the proof.
\end{proof}
\begin{theorem} \label{k_theorem}
Using the notation of Lemma~\ref{kolmo-inequality}, assume that 
\beq \no
\sum_{n=1}^\infty \mean{z_\omega(n)^2} <\infty.
\eeq
Then
\beq \no
\sum_{n=1}^\infty z_\omega(n)
\eeq
converges almost surely. Furthermore, for any $n$,
\beq \label{bound_on_sum}
\mean{ \left( \sum_{j=n}^\infty z_\omega(j) \right)^2}\leq \sum_{j=1}^\infty \mean{z_\omega(j)^2} <\infty.
\eeq
\end{theorem}
\begin{proof}
For any $\veps>0$ and for any $N_1<N_2$, by \eqref{k-inequality},
\beq \no
P \left( \left \{\omega \;\bigg|\; \max_{N_1\leq n \leq N_2} |z_\omega(n)+\ldots+z_\omega(N_2)|> \veps \right \} \right)\leq
\frac{\sum_{j=N_1}^{N_2} \mean{ z_\omega(j)^2}}{\veps^2}.
\eeq
Thus, the event 
\beq \no
\left \{ \omega \;\bigg|\; \exists N_1,N_2, \textrm{ arbitrarily large,}
\max_{N_1\leq n \leq N_2} |z_\omega(n)+\ldots+z_\omega(N_2)|> \veps \right \}
\eeq
has probability zero.
So we get that, with probability 
one, for any $\veps>0$, there exists $N_{\veps}$ so that for any $N_1>N_2>
N_{\veps}$, 
\beq \no
\left|\sum_{j=N_1}^{N_2}z_\omega(j)\right|<\veps, 
\eeq
or, in other words, 
$\sum_{j=1}^\infty z_\omega(j)$ converges with probability one. \eqref{bound_on_sum} now follows from 
Fatou's lemma. 
\end{proof}

\section{A Central Lemma}

The idea at the basis of our analysis is that of variation of 
parameters. We want to obtain a `linear' relationship between the generalized
eigenfunctions of the original problem and those of the perturbed problem. 
Thus, for fixed $E\in\bbR$, let $T^E_0(n)$ and $S^E_0(n)$ denote the 
$n$-steps and one-step transfer matrices respectively and let 
$T^E_\omega(n)$ and $S^E_\omega(n)$ denote the same objects for the 
perturbed problem (depending on the random parameter $\omega$). Define 
$\calD^E_\omega(n)$ through the equation:
\beq \label{Anomega-definition}
T^E_\omega(n)=T^E_0(n) \calD^E_\omega(n).
\eeq
Then
\beq \label{Aomega-recursion}
\begin{split}
\calD^E_\omega(n-1) &=T^E_0(n-1)^{-1}T^E_\omega(n-1) \\
&=T^E_0(n-1)^{-1}T^E_\omega(n-1)T^E_\omega(n)^{-1} T^E_0(n)T^E_0(n)^{-1}T^E_\omega(n)\\
&=T^E_0(n-1)^{-1}T^E_\omega(n-1)T^E_\omega(n)^{-1}T^E_0(n)\calD^E_\omega(n) \\
&=T^E_0(n)^{-1}S^E_0(n)S^E_\omega(n)^{-1}T^E_\omega(n) T^E_\omega(n)^{-1}T^E_0(n)\calD^E_\omega(n)\\
&=T^E_0(n)^{-1}S^E_0(n)S^E_\omega(n)^{-1}T^E_0(n)\calD^E_\omega(n)\\
&=(I+\calU^E_\omega(n))\calD^E_\omega(n),
\end{split}
\eeq
where
\beq \label{Unomega-definition}
\calU^E_\omega(n)=T^E_0(n)^{-1}(S^E_0(n)S^E_\omega(n)^{-1}-I)T^E_0(n).
\eeq
Almost sure convergence of $\calD^E_\omega(n)$ to $I$ would insure that the 
asymptotic properties of $T^E_\omega(n)$ would resemble 
those of $T^E_0(n)$. This would suffice in the absolutely continuous case. 
In the singular continuous case we would like to control the convergence 
rate of each of the column vectors of $\calD^E_\omega(n)$ separately. The 
following lemma is actually a random version of a well known result on the 
control of the amplitudes (see for instance \cite{kls} and problem XI.97 
in \cite{reed-simon3}). It is central to everything that follows.

\begin{lemma} \label{random_variation}
Let $\left \{\calU(n)=\left( \begin{array}{cc}
u_{11}(n) & u_{12}(n) \\
u_{21}(n) & u_{22}(n)
\end{array} \right) \right \}_{n=1}^\infty$ be a sequence of matrices in
$M_2(\bbR)$, and let
$\se{\ti{b}_\omega}$ be a sequence of independent random variables with 
zero mean. Suppose that
\beq \label{decay_condition}
\sum_{n=1}^\infty
\left( \mean{\ti{b}_\omega(n)^2} (u_{11}(n)^2+u_{12}(n)^2+u_{22}(n)^2+u_{21}(n)^2f_+(n)^2) \right) 
<\infty 
\eeq
for some monotonically increasing sequence - $f_+(n)>0$.
Then, $P$-almost surely,
\beq \label{conclusion_for_a}
\bfd_\omega(n-1)-\bfd_\omega(n)=\ti{b}_\omega(n)\calU(n)\bfd_\omega(n)
\eeq
has solutions - 
$\bfd^+_\omega(n) \equiv \left( \begin{array} {c} d^+_{1,\omega}(n) \\
d^+_{2,\omega}(n)
\end{array} \right)$ 
and 
$\bfd^-_\omega(n) \equiv \left ( \begin{array} {c} d^-_{1,\omega}(n) \\
d^-_{2,\omega}(n)
\end{array} \right)$, 
that satisfy
\begin{eqnarray}
\label{asymptotic_amplitude1}
\lim_{n \rightarrow \infty} d_{1,\omega}^+(n) & = & 0 \\
\label{asymptotic_amplitude2}
\lim_{n \rightarrow \infty} d_{2,\omega}^+(n) & = & 1 \\
\label{asymptotic_amplitude3}
\lim_{n \rightarrow \infty} d_{1,\omega}^-(n) & =  & 1 \\
\label{asymptotic_amplitude4}
\lim_{n \rightarrow \infty} d_{2,\omega}^-(n)f_+(n) & =  & 0
\end{eqnarray}
\end{lemma}
\begin{remark}
Note that, for $f_+(n) \equiv 1$, \eqref{asymptotic_amplitude1}-\eqref{asymptotic_amplitude4} mean that the matrix 
equation 
\beq \no
\calD_\omega(n-1)= \left(I+\ti{b}_\omega(n) \calU(n) \right) \calD_\omega(n)
\eeq
has a solution $\calD_\omega(n)$ such that $\lim_{n\rightarrow \infty} \calD_\omega(n)=I$. 
\end{remark}
\begin{proof}
We start by constructing $\bfd^+$. Let
\beq \label{a+0}
\bfd^{+,0}_\omega(n) \equiv \left( \begin{array}{c}
0 \\
1 \end{array} \right)
\eeq
and denote
\beq \label{Utilde}
\ti{\calU}_\omega(n)=\ti{b}_\omega(n) \calU(n).
\eeq
Then Theorem~\ref{k_theorem} says that
\beq \label{a+1}
\bfd^{+,1}_\omega(n)=\sum_{j=n+1}^\infty \ti{\calU}_\omega(j)\bfd^{+,0}_\omega(j)
\eeq
is defined $P$-a.s.\ for any $n$ and that 
$\mean{\parallel \bfd^{+,1}_\omega(n) \parallel ^2} $
is bounded in $n$. Note also, that $\bfd^{+,1}_\omega(n)$ is a measurable
function of
\beq \no
\{\ti{b}_\omega(n+1),\ti{b}_\omega(n+2),...\}.
\eeq
Now, for $k \geq 1$, assume that
$\bfd^{+,k}_\omega(n)$ is defined $P$-a.s.\ as a measurable function of
\beq \no
\{\ti{b}_\omega(n+1),\ti{b}_\omega(n+2),...\}
\eeq
and that
\mbox{$ \mean{\parallel \bfd^{+,k}_\omega(n) \parallel ^2} $} is bounded in $n$.
Then by Theorem~\ref{k_theorem}, it is possible to define
\beq \label{a+(k+1)}
\bfd^{+,(k+1)}_\omega(n)=\sum_{j=n+1}^\infty
\ti{\calU}_\omega(j)\bfd^{+,k}_\omega(j)
\eeq
$P$-a.s.\ and this definition satisfies all of the
properties listed above. Thus, by induction, we construct
 $\bfd^{+,k}_\omega(n)$ for every $k\in\bbN$.
Now,
\begin{align} \no
&\mean{ \parallel \bfd^{+,k}_\omega(n) \parallel ^2 } \notag \\
&=\mean{\parallel  \sum_{j=n+1}^\infty \ti{\calU}_\omega(j)(j)\bfd^{+,(k-1)}_\omega(j) \parallel ^2} \notag \\
& =\mean{\left(\sum_{j=n+1}^\infty \ti{b}_\omega(j)\left( 
u_{11}(j)d^{+,(k-1)}_{1,\omega}(j)+ u_{12}(j)d^{+,(k-1)}_{2,\omega}(j)
\right) \right)^2 } \notag \\
& \quad +\mean{\left(\sum_{j=n+1}^\infty \ti{b}_\omega(j)\left( 
u_{21}(j)d^{+,(k-1)}_{1,\omega}(j) + u_{22}(j)d^{+,(k-1)}_{2,\omega}(j)
\right)   \right)^2 } \notag \\
& \leq \sum_{j=n+1}^\infty \mean{\ti{b}_\omega(j)^2 } \mean{\left( 
u_{11}(j)d^{+,(k-1)}_{1,\omega}(j) + u_{12}(j)d^{+,(k-1)}_{2,\omega}(j)
\right)^2 } \notag \\
& \quad + \sum_{j=n+1}^\infty \mean{ \ti{b}_\omega(j)^2 } \mean{\left( 
u_{21}(j)d^{+,(k-1)}_{1,\omega}(j) + u_{22}(j)d^{+,(k-1)}_{2,\omega}(j)
\right)^2 } \notag \\
& = \sum_{j=n+1}^\infty \mean{ \left \| \ti{\calU}_\omega(j)\bfd^{+,(k-1)}_\omega(j) \right \| ^2 }, \notag
\end{align}
where, for the inequality, we used Fatou's lemma, the independence of the $\ti{b}_\omega(j)$ and the fact that 
$\bfd^{+,(k-1)}_\omega(j)$ is a function of $\{ \ti{b}_\omega(j+1), \ti{b}_\omega(j+2) \ldots \}$ only. 
Now, there exists a
universal constant $C$, such that for any $2 \times 2$ matrix $\mathcal{A}$,
\beq \no
\parallel \mathcal{A} \parallel ^2 \leq C \parallel \mathcal{A} \parallel ^2_{\rm HS},
\eeq 
where $ \parallel \cdot \parallel _{\rm HS}$ is the Hilbert-Schmidt norm. Therefore, using independence again,
\begin{align} \no
&\mean{ \left \| \ti{\calU}_\omega(j)\bfd^{+,(k-1)}_\omega(j) \right \| ^2 } \leq 
\mean{ \left \| \ti{\calU}_\omega(j) \right \| ^2 \left \| \bfd^{+,(k-1)}_\omega(j) \right \| ^2 } \notag \\ 
& \leq C \mean{\left \| \ti{\calU}_\omega(j) \right \| ^2_{\rm HS} \left \| \bfd^{+,(k-1)}_\omega(j) \right \| ^2} \notag \\
& =C \mean{\left\| \ti{\calU}_\omega(j) \right \| ^2_{\rm HS}} 
\mean{\left \| \bfd^{+,(k-1)}_\omega(j) \right \| ^2}. \notag 
\end{align}
Thus,
\beq \no 
\mean{\left\| \bfd^{+,k}_\omega(n) \right\| ^2 } \leq \left( \sup_{j>n}\mean{\left \| \bfd^{+,(k-1)}_\omega(j) \right\| ^2 } \right)
\cdot C \sum_{j=n+1}^\infty \mean{\left\| \ti{\calU}_\omega(j) \right\| ^2_{\rm HS}}, 
\eeq
and therefore, for any $N\in\bbN,$
\beq \no
\begin{split}
\sup_{n \geq N} \mean{\left\| \bfd^{+,k}_\omega(n) \right\| ^2 } & \leq \left( \sup_{n \geq 
N} \mean{\left\| \bfd^{+,(k-1)}_\omega(n) \right\| ^2 } \right) \\ 
& \times C\sum_{j=N+1}^\infty \left( \mean{\left\| \ti{\calU}_\omega(j) \right\|^2_{\rm HS}} \right)
\end{split}
\eeq
with $C$ independent of $k$ and $N$. Thus, it follows that 
\beq \no
\mean{ \left\| \bfd^{+,k}_\omega(N) \right\| ^2 } \leq \left( C\sum_{j=N+1}^\infty  
\mean{(\ti{b}_\omega(j))^2} \left \| \calU(j) \right\|^2_{\rm HS} \right)^k.
\eeq

 From \eqref{decay_condition}, it thus follows that there exists some 
$N\in\bbN$,
which we denote by $N_{1/4}$, so that for any $n \geq N_{1/4}$
\beq \label{N1/4_condition}
\mean{\left\| \bfd^{+,k}_\omega(n) \right\| ^2} \leq \left(\frac{1}{4}\right)^k.
\eeq
Let
\beq \no
\Omega_{k,n}=\left \{\omega \left \vert \left\| \bfd^{+,k}_\omega(n) \right\|  \geq \left(\frac{1}{4}\right)^{\frac{k}{4}} \right. 
\right \}.
\eeq
Then, by Chebyshev's inequality, for $n \geq N_{1/4}$
\beq \no
P(\Omega_{k,n}) \leq \left(\frac{1}{2}\right)^k.
\eeq
Thus, by the Borel Cantelli lemma, for $n$ large enough, there exists a 
set $\Omega^0 \subseteq \Omega$ of full $P$ measure such that for any 
$\omega\in \Omega^0$
\beq \label{a+definition}
\bfd^+_\omega(n)=\sum_{k=0}^{\infty}\bfd^{+,k}_\omega(n)
\eeq
is defined.

Suppose for a while that we could show
\beq \label{convergence_and_equality}
\sum_{k=0}^\infty \sum_{j=n+1}^\infty \ti{\calU}_\omega(j)\bfd^{+,k}_\omega(j)=
\sum_{j=n+1}^\infty \sum_{k=0}^\infty \ti{\calU}_\omega(j)\bfd^{+,k}_\omega(j)
\eeq
$P$-a.s.\ and for large enough $n$, in the 
sense 
that 
both sides converge and are equal. Then we would have
\beq \no
\begin{split}
\bfd^+_\omega(n)-\left(\begin{array}{c}
0 \\
1
\end{array} \right) & =\sum_{k=1}^\infty \bfd^{+,k}_\omega(n)=\sum_{k=0}^\infty \sum_{j=n+1}^\infty \ti{\calU}_\omega(j)\bfd^{+,k}_\omega(j) \\
&=\sum_{j=n+1}^\infty \sum_{k=0}^\infty \ti{\calU}_\omega(j)\bfd^{+,k}_\omega(j) 
= \sum_{j=n+1}^\infty \ti{\calU}_\omega(j) \bfd^+_\omega(j),
\end{split}
\eeq
which implies
\beq \no
\bfd^+_\omega(n-1)-\bfd^+_\omega(n)=\ti{\calU}_\omega(n)\bfd^+_\omega(n),
\eeq
which is \eqref{conclusion_for_a}. Furthermore, 
\eqref{asymptotic_amplitude1} and \eqref{asymptotic_amplitude2} would be 
obvious from the convergence.

Therefore, we need to prove \eqref{convergence_and_equality}. We know that for $n \geq N_{\frac{1}{4}}$ 
\beq \no
\sum_{k=0}^\infty \bfd^{+,k}_\omega(n)=\sum_{k=0}^\infty \sum_{j=n+1}^\infty \ti{\calU}_\omega(j) \bfd^{+,k}_\omega(j)
\eeq 
converges $P$-a.s.\ which is precisely the convergence of the LHS. It is also obvious 
that 
$\bfd^+_\omega(n)$ is a measurable function of 
$\{\ti{b}_\omega(n+1),\ti{b}_\omega(n+2)...\}$, so if we show uniform 
boundedness 
of $\mean{\parallel \bfd^+_\omega(n) \parallel ^2} $, we will have the convergence of 
the RHS$=\sum_{j=n+1}^\infty \ti{\calU}_\omega(j) \bfd^+_\omega(j)$, by 
Theorem~\ref{k_theorem}. But this is true, since
\beq \no
\begin{split}
\mean{\left\| \bfd^+_\omega(n) \right\| ^2} &= \mean{\lim_{N \rightarrow \infty} 
\left\| \left( \sum_{k=0}^N \bfd^{+,k}_\omega(n) \right) \right\| ^2 } \\
&\leq \liminf_{N \rightarrow \infty}
\mean{\left\| \left( \sum_{k=0}^N \bfd^{+,k}_\omega(n) \right) \right\| ^2} \\
&\leq \liminf_{N \rightarrow \infty} 
\mean{\left( \sum_{k=0}^N \left \| \bfd^{+,k}_\omega(n) \right\| \right)^2}  \\
& = \liminf_{N \rightarrow \infty} 
\Big \langle \left( \left\| \bfd^{+,0}_\omega(n) \right\| + \ldots +
\left\| \bfd^{+,N}_\omega(n) \right\| \right) \\
& \quad \times \left( \left\| \bfd^{+,0}_\omega(n) \right\| + \ldots + 
\left\| \bfd^{+,N}_\omega(n) \right\| \right) \Big \rangle \\
& = \liminf_{N \rightarrow \infty} \Big( \mean{ \left\| \bfd^{+,0}_\omega(n) \right\| ^2 }+
2 \mean{ \left\| \bfd^{+,1}_\omega(n) \right\| \left\| \bfd^{+,0}_\omega(n) \right\| }  \\
& \quad +\mean{ \left\| \bfd^{+,1}_\omega(n) \right\| ^2} + \ldots \Big) 
\end{split}
\eeq
where, in the fourth and fifth lines, each factor of the form $\|\bfd^{+,k_0}\|$ in one set of summands, is coupled
to factors of the form $\|\bfd^{+,k}\|$ for $k \leq k_0$ in the other set of summands. 
Using this way of writing the product, the Cauchy-Schwarz inequality,
and the fact that, for $n \geq N_{\frac{1}{4}}$, \eqref{N1/4_condition} holds (so that, in particular, all factors are 
bounded by 1 from above), we get that
\beq \label{uniform_boundedness}
\mean{ \parallel \bfd^+_\omega(n) \parallel ^2 } \leq 
\liminf_{N \rightarrow \infty} \sum_{k=0}^N (2k+1)\left(\frac{1}{2}\right)^k,
\eeq
so that $\mean{ \parallel \bfd^+_\omega(n) \parallel ^2 } $ is bounded in $n$.
Thus we are left with proving the equality 
\eqref{convergence_and_equality}, or in other words, with proving
\beq \no
\lim_{K\rightarrow\infty}\left( \sum_{j=n+1}^\infty\sum_{k=0}^\infty\ti{\calU}_\omega(j)
\bfd^{+,k}_\omega(j)-\sum_{j=n+1}^\infty\sum_{k=0}^K\ti{\calU}_\omega(j)
\bfd^{+,k}_\omega(j) \right)=0
\eeq
which is the same as
\beq \label{convergence_to_zero} 
\lim_{K\rightarrow\infty}\sum_{j=n+1}^\infty\sum_{k=K+1}^\infty\ti{\calU}_\omega(j)
\bfd^{+,k}_\omega(j)=0
\eeq
$P$-almost surely. Denote
\beq \label{akj}
\mathfrak{d}^{+,K}_\omega(j)=\sum_{k=K+1}^\infty \bfd^{+,k}_\omega(j).
\eeq
Then,
\beq \no      
\begin{split}
&P \left \{ \left \| \sum_{j=n+1}^\infty \ti{\calU}_\omega(j)\mathfrak{d}^{+,K}_\omega(j) \right \|  \geq 
\frac{1}{K} \right \} \\
&\leq K^2 \mean{\left \| \sum_{j=n+1}^\infty 
\ti{\calU}_\omega(j)\mathfrak{d}^{+,K}_\omega(j) \right \| ^2 } \\
&\leq K^2 \liminf_{N\rightarrow\infty} 
\sum_{j=n+1}^N \mean{\left \| \ti{\calU}_\omega(j)\mathfrak{d}^{+,K}_\omega(j) \right \| ^2 } \\
&\leq K^2 \left( \sum_{j=n+1}^\infty \mean{\left \| \ti{\calU}_\omega(j) \right \| ^2 }\right) 
\sup_{j \geq n+1} 
\mean{\left \| \mathfrak{d}^{+,K}_\omega(j)\right \| ^2}=**. 
\end{split}
\eeq
The same considerations that lead to \eqref{uniform_boundedness}, lead to the 
conclusion that $**\in \ell^1(K)$ and therefore, by Borel Cantelli 
\eqref{convergence_to_zero} holds, $P$-almost surely, and we are done.
\ \\

To construct $\bfd^-$ go through the same procedure, constructing $\bfd^{-,k}_\omega(n)$, with $\bfd^{-,0} \equiv 
\left( \begin{array}{c} 
1 \\
0
\end{array}
\right)$, and define
\beq \no
\bfd^-_\omega(n)=\sum_{k=0}^\infty \bfd^{-,k}_\omega(n)
\eeq
for $n$ large enough.
Everything works the same as for the construction of $\bfd^+_\omega(n)$.

To show \eqref{asymptotic_amplitude4}, we define 
\beq \label{Xn}
\calX(n) \equiv \left( \begin{array}{cc}
1 & 0 \\
0 & f_+(n)
\end{array} \right)
\eeq
and note that for $m>n$,
\beq \label{monotonicity}
\left\| \calX(n)\calX(m)^{-1}\right\| \leq 1.
\eeq
Denote now
\beq \label{Wjomega}
\calW_\omega(n)=\calX(n)\ti{\calU}_\omega(n)\calX(n)^{-1},
\eeq
and 
\beq \no
\ti{\bfd}^{-,k}_\omega(n)=\calX(n)\bfd^{-,k}_\omega(n),
\eeq
so that we have 
\begin{align} \no
& \ti{\bfd}^{-,k}_\omega(n) \notag \\
& =\calX(n) \sum_{j=n+1}^\infty \ti{\calU}_\omega(j)\bfd^{-,(k-1)}_\omega(j) \notag \\
& =\calX(n) \sum_{j=n+1}^\infty \ti{\calU}_\omega(j)\calX(j)^{-1}\calX(j)\bfd^{-,(k-1)}_\omega(j) \notag \\
& =\calX(n) \sum_{j=n+1}^\infty \calX(j)^{-1}\calX(j)\ti{\calU}_\omega(j)\calX(j)^{-1} \ti{\bfd}^{-,(k-1)}_\omega(j) \notag \\
& =\sum_{j=n+1}^\infty \calX(n)\calX(j)^{-1}\calW_\omega(j) \ti{\bfd}^{-,(k-1)}_\omega(j). \notag
\end{align}
Since, by \eqref{decay_condition},
\beq \label{decay2}
\sum_{n=1}^{\infty} \mean{ \left \| \calW_\omega(n) \right \| ^2 } <\infty,
\eeq
and by \eqref{monotonicity}
\beq \label{uniform-matrix-convergence}
\sum_{j=n+1}^{\infty} \mean{ \left \| \calX(n)\calX(j)^{-1}\calW_\omega(j) \right \| ^2 }
\leq \sum_{j=n+1}^{\infty} \mean{ \left \| \calW_\omega(j) \right \| ^2 },
\eeq
we can repeat the argument in the first part of the proof to show that
\beq \label{uniform-boundedness-for-tildea}
\sup_j \mean{ \left \| \ti{\bfd}^-_\omega(j) \right \| ^2 } 
\equiv \sup_j \mean{\left \| \calX(j)\bfd^-_\omega(j) \right \| ^2}  <\infty. 
\eeq

Now,
\beq \label{X(n)a-estimate}
\begin{split}
\left\| \calX(n)\bfd^-_\omega(n)-\left( \begin{array}{c}
1 \\
0
\end{array} \right)\right\| &=
\left \| \calX(n)\bfd^-_\omega(n)-\calX(n)\left( \begin{array}{c}
1 \\
0
\end{array} \right)\right \| \\
&= \left \| \calX(n)\sum_{j=n+1}^\infty 
\ti{\calU}_\omega(j)\bfd^-_\omega(j) \right \| \\
& =\left \| \sum_{j=n+1}^\infty \calX(n)\calX(j)^{-1}\calW_\omega(j)\calX(j)\bfd^-_\omega(j)\right \|.
\end{split}
\eeq

Since $ \calX(j)\bfd^-_\omega(j) $ is a measurable function of $\{\ti{b}_\omega(j+1),\ti{b}_\omega(j+2),\ldots\}$, and, 
since
$\calX(n)\calX(j)^{-1}\calW_\omega(j)=\ti{b}_\omega(j) \hat{\calW}(n,j)$ (where $\hat{\calW}(n,j)$ is a deterministic 
matrix), 
one can repeat the
proof of Lemma~\ref{kolmo-inequality}, using \eqref{uniform-matrix-convergence} and \eqref{uniform-boundedness-for-tildea},
to show that  
\beq \no
\lim_{n\rightarrow\infty}
\left \| \sum_{j=n+1}^\infty \calX(n)\calX(j)^{-1}\calW_\omega(j)\calX(j)\bfd^-_\omega(j) \right \| =0
\eeq
$P$-almost surely, which, by \eqref{X(n)a-estimate}, is exactly \eqref{asymptotic_amplitude3} and 
\eqref{asymptotic_amplitude4}. 
\end{proof}

As remarked earlier, this lemma is actually a `random variation' on a deterministic stability result. This random version
uses the zero mean of the random variables in order to replace an $\ell^1$ summability condition (which is the natural 
condition in the deterministic case) with an $\ell^2$ condition. It is natural to ask whether it is possible to obtain 
such a result for a situation in which there is a combination of terms -- coefficients which are $\ell^2$ with zero mean 
and coefficients that are $\ell^1$. The following lemma is an extension of Lemma~\ref{random_variation} in this direction
(in the special case `$f_+ \equiv 1$'), which is tailored especially for our needs in the next section.

\begin{lemma} \label{random_variation2}
Let 
\beq \no
\left \{\calU(n)=\left( \begin{array}{cc}
u_{11}(n) & u_{12}(n) \\
u_{21}(n) & u_{22}(n)
\end{array} \right) \right \}_{n=1}^\infty,
\eeq
\beq \no
\left \{\calV(n)=\left( \begin{array}{cc}
v_{11}(n) & v_{12}(n) \\
v_{21}(n) & v_{22}(n)
\end{array} \right) \right \}_{n=1}^\infty,
\eeq
\beq \no
\left \{\calW(n)=\left( \begin{array}{cc}
w_{11}(n) & w_{12}(n) \\
w_{21}(n) & w_{22}(n)
\end{array} \right) \right \}_{n=1}^\infty
\eeq
be three matrix-valued sequences, and let
$\se{\ti{b}_{1,\omega}}$, $\se{\ti{b}_{2,\omega}}$, $\se{\ti{b}_{3,\omega}}$ 
be three sequences of random variables that satisfy the following properties:
\begin{enumerate}
\item For any $i,j=1,2,3$ and $n_1 \neq n_2$, $\ti{b}_{i,\omega}(n_1)$ and $\ti{b}_{j,\omega}(n_2)$ are independent random 
variables. 
\item For any $n$, 
\beq \label{mean-zero}
\mean{\ti{b}_{1,\omega}(n)} = \mean{\ti{b}_{2,\omega}(n)}= 
\mean{\ti{b}_{1,\omega}(n) \ti{b}_{2,\omega}(n)} = \mean{\ti{b}_{1,\omega}(n) \ti{b}_{3,\omega}(n)} =0 
\eeq
\item For all $n$ and any $\omega \in \Omega$, $\ti{b}_{3,\omega}(n) \geq 0$.
\item 
\beq \label{decay-condition-b1}
\sum_{n=1}^\infty \mean{\ti{b}_{1,\omega}(n)^2} \Vert \calU(n) \Vert ^2_{\rm HS} < \infty
\eeq
and
\beq \label{decay-condition-b2}
\sum_{n=1}^\infty \mean{\ti{b}_{2,\omega}(n)^2} \Vert \calV(n)\Vert ^2_{\rm HS} < \infty.
\eeq
\item 
\beq \label{decay-condition-b3}
\sum_{n=1}^\infty \mean{\ti{b}_{3,\omega}(n)^2} ^{1/2} \Vert \calW(n)\Vert _{\rm HS} < \infty.
\eeq
\item
\beq \label{decay-condition-b2b3}
\sum_{n=1}^\infty \mean{ |\ti{b}_{2,\omega}(n)| \ti{b}_{3,\omega}(n)} \Vert \calV(n)\Vert _{\rm HS} 
\Vert \calW(n)\Vert _{\rm HS} <\infty.
\eeq
\end{enumerate}

Then, $P$-almost surely,
\beq \label{conclusion_for_a1}
\bfd_\omega(n)-\bfd_\omega(n-1)=\left( \ti{b}_{1,\omega}(n)\calU(n)+\ti{b}_{2,\omega}(n) \calV(n)+
\ti{b}_{3,\omega}(n) \calW(n) \right)\bfd_\omega(n)
\eeq
has solutions - 
$\bfd^+_\omega(n) \equiv \left ( \begin{array} {c} d^+_{1,\omega}(n) \\
d^+_{2,\omega}(n)
\end{array} \right)$ 
and 
$\bfd^-_\omega(n) \equiv \left ( \begin{array} {c} d^-_{1,\omega}(n) \\
d^-_{2,\omega}(n)
\end{array} \right)$, 
that satisfy
\begin{eqnarray}
\label{asymptotic_amplitude12}
\lim_{n \rightarrow \infty} d_{1,\omega}^+(n) & = & 0 \\
\label{asymptotic_amplitude22}
\lim_{n \rightarrow \infty} d_{2,\omega}^+(n) & = & 1 \\
\label{asymptotic_amplitude32}
\lim_{n \rightarrow \infty} d_{1,\omega}^-(n) & =  & 1 \\
\label{asymptotic_amplitude42}
\lim_{n \rightarrow \infty} d_{2,\omega}^-(n) & =  & 0
\end{eqnarray}
\end{lemma}

\begin{proof}
This proof follows the same strategy of the proof of Lemma~\ref{random_variation}. We shall try to avoid unnecessary 
repetitions. The first step is the construction of $\bfd^{+,k}$ for $k \geq 0$. As before, let
\beq \label{a+0,2}
\bfd^{+,0}_\omega(n) \equiv \left( \begin{array}{c}
0 \\
1 \end{array} \right).
\eeq 
Now, note that, by H\"older,
\beq \label{b3-absolute-convergence1}
\mean{\ti{b}_{3,\omega}(n)} \leq \mean{\ti{b}_{3,\omega}(n)^2 }^{1/2}, 
\eeq
so
\beq \label{b3-absolute-convergence2}
\sum_{n=1}^\infty \mean{ \ti{b}_{3,\omega}(n) } \Vert \calW(n)\Vert _{\rm HS} 
\leq \sum_{n=1}^\infty \mean{\ti{b}_{3,\omega}(n)^2}^{1/2} \Vert \calW(n)\Vert _{\rm HS}<\infty. 
\eeq 
As in the preceding proof, we want to show that if 
$\bfd^{+,k}_\omega(n)$ is defined $P$-a.s.\ as a measurable function of $\{\ti{b}_{i,\omega}(j)\}_{i=1,2,3,\ j>n}$, 
for any $n$; and $\mean{\Vert \bfd^{+,k}_\omega(n)\Vert^2}$ is bounded in $n$, 
then the same holds true for 
\beq \label{a+(k+1),2}
\bfd^{+,(k+1)}_\omega(n)=\sum_{j=n+1}^\infty \left( \ti{b}_{1,\omega}(j)\calU(j)+\ti{b}_{2,\omega}(j)\calV(j)+
\ti{b}_{3,\omega}(j)\calW(j) 
\right) \bfd^{+,k}_\omega(j).
\eeq
If indeed $\mean{ \Vert \bfd^{+,k}_\omega(n)\Vert ^2 }$ is bounded in $n$ (say, by $C$), then, by H\"older, so is 
$\mean{\Vert \bfd^{+,k}_\omega (n)\Vert}$ so by \eqref{b3-absolute-convergence2} and from the independence we will have
\beq \label{b3-absolute-convergence4}
\begin{split}
\sum_{n=1}^\infty \mean{\left\| \ti{b}_{3,\omega}(n) \calW(n) \bfd^{+,k}_\omega(n)\right \|} &\leq
\sum_{n=1}^\infty \mean{\ti{b}_{3,\omega}(n) \left\| \calW(n)\right\|  \left\| \bfd^{+,k}_\omega(n)\right\|} \\
&=\sum_{n=1}^\infty \mean{\ti{b}_{3,\omega}(n) } \left\| \calW(n)\right\|  \mean{\left\| \bfd^{+,k}_\omega(n)\right\|} \\
& \leq C \sum_{n=1}^\infty \mean{\ti{b}_{3,\omega}(n)} \left\| \calW(n)\right\|  < \infty.     
\end{split}
\eeq
Therefore, monotone convergence implies that 
\beq \label{b3-absolute-convergence5}
\sum_{n=1}^\infty \ti{b}_{3,\omega}(n) \calW(n) \bfd^{+,k}_\omega(n)
\eeq
is absolutely convergent $P$-a.s. Theorem~\ref{k_theorem} implies the almost sure convergence of the first two 
summands in \eqref{a+(k+1),2} so $\bfd^{+,(k+1)}_\omega(n)$ is defined $P$-a.s.\ as a measurable function of 
$\{\ti{b}_{i,\omega}(j)\}_{i=1,2,3,\ j>n}$. Thus, we are left with showing that $\mean{\Vert \bfd^{+,(k+1)}_\omega(n)\Vert ^2}$ 
is bounded in $n$. We proceed to estimate
\begin{align} \label{a+(k+1)^2,2basic}
&\mean{ \left\| \bfd^{+,(k+1)}_\omega(n)\right\| ^2} \notag \\
& \leq 
\liminf_{N\rightarrow \infty} \Big \langle \Big \|  \sum_{j=n+1}^N \ti{b}_{1,\omega}(j)\calU(j)\bfd^{+,k}_\omega(j) \notag \\
&\quad +\sum_{j=n+1}^N \ti{b}_{2,\omega}(j)\calV(j)\bfd^{+,k}_\omega(j)+
\sum_{j=n+1}^N \ti{b}_{3,\omega}(j)\calW(j)\bfd^{+,k}_\omega(j)\Big\| ^2 \Big \rangle 
\notag \\
&=\liminf_{N \rightarrow \infty} \Bigg( \bigg \langle \Big |\sum_{j=n+1}^N 
\Big( \ti{b}_{1,\omega}(j)(u_{11}(j)d^{+,k}_{1,\omega}(j)+u_{12}(j)d^{+,k}_{2,\omega}(j)) \notag \\
& \quad+ \ti{b}_{2,\omega}(j)(v_{11}(j)d^{+,k}_{1,\omega}(j)+v_{12}(j)d^{+,k}_{2,\omega}(j)) \notag \\
& \quad+\ti{b}_{3,\omega}(j)(w_{11}(j)d^{+,k}_{1,\omega}(j)+w_{12}(j)d^{+,k}_{2,\omega}(j))\Big) \Big|^2 \bigg \rangle \notag \\
& \quad+ \bigg \langle \Big| \sum_{j=n+1}^N  \Big( 
\ti{b}_{1,\omega}(j)(u_{21}(j)d^{+,k}_{1,\omega}(j)+u_{22}(j)d^{+,k}_{2,\omega}(j)) 
\notag \\
&\quad + \ti{b}_{2,\omega}(j)(v_{21}(j)d^{+,k}_{1,\omega}(j)+v_{22}(j)d^{+,k}_{2,\omega}(j)) \notag \\
& \quad +\ti{b}_{3,\omega}(j)(w_{21}(j)d^{+,k}_{1,\omega}(j)+w_{22}(j)d^{+,k}_{2,\omega}(j)) \Big) \Big|^2 \bigg \rangle \Bigg) 
\notag \\
&\equiv \liminf_{N \rightarrow \infty} \Big( R_1(n,N)+R_2(n,N) \Big). 
\end{align}
Now, using independence and \eqref{mean-zero} we see that
\begin{align} \label{a+(k+1)^2,2stagetwo}
&R_1(n,N) \notag \\
&=\bigg \langle \Big |\sum_{j=n+1}^N 
\Big( \ti{b}_{1,\omega}(j)(u_{11}(j)d^{+,k}_{1,\omega}(j)+u_{12}(j)d^{+,k}_{2,\omega}(j)) \notag \\
&\quad + \ti{b}_{2,\omega}(j)(v_{11}(j)d^{+,k}_{1,\omega}(j)+v_{12}(j)d^{+,k}_{2,\omega}(j)) \notag \\
& \quad +\ti{b}_{3,\omega}(j)(w_{11}(j)d^{+,k}_{1,\omega}(j)+w_{12}(j)d^{+,k}_{2,\omega}(j))\Big) \Big|^2 \bigg \rangle \notag \\
&=\sum_{j=n+1}^N \mean{\ti{b}_{1,\omega}(j)^2} 
\mean{\left(u_{11}(j)d^{+,k}_{1,\omega}(j)+u_{12}(j)d^{+,k}_{2,\omega}(j)\right)^2} \notag \\
& \quad +\sum_{j=n+1}^N \mean{\ti{b}_{2,\omega}(j)^2} 
\mean{\left(v_{11}(j)d^{+,k}_{1,\omega}(j)+v_{12}(j)d^{+,k}_{2,\omega}(j)\right)^2 } \notag \\
& \quad +\sum_{j=n+1}^N \Big(\mean{\ti{b}_{2,\omega}(j) \ti{b}_{3,\omega}(j) } \notag \\  
& \quad \times \mean{ (v_{11}(j)d^{+,k}_{1,\omega}(j)+v_{12}(j)d^{+,k}_{2,\omega}(j))
(w_{11}(j)d^{+,k}_{1,\omega}(j)+w_{12}(j)d^{+,k}_{2,\omega}(j)) } \Big) \notag \\
& \quad +\mean{ \left( \sum_{j=n+1}^N \left( \ti{b}_{3,\omega}(j)
\left(w_{11}(j)d^{+,k}_{1,\omega}(j)+w_{12}(j)d^{+,k}_{2,\omega}(j)\right) \right) \right)^2 } \notag \\
& \equiv R_1^I(n,N)+R_1^{II}(n,N)+R_1^{III}(n,N)+R_1^{IV}(n,N)
\end{align}
and a similar expression holds for $R_2(n,N)$. Thus we have:
\beq \no
\begin{split}
&R_1(n,N)+R_2(n,N)=\underbrace{R_1^I(n,N)+R_2^I(n,N)}_{R^I(n,N)} +
\underbrace{R_1^{II}(n,N)+R_2^{II}(n,N)}_{R^{II}(n,N)} \\
&+\underbrace{R_1^{III}(n,N)+R_2^{III}(n,N)}_{R^{III}(n,N)} +
\underbrace{R_1^{IV}(n,N)+R_2^{IV}(n,N)}_{R^{IV}(n,N)}.
\end{split}
\eeq
\eqref{decay-condition-b1} means that
\beq \label{bound-on-RI}
\begin{split}
R^I(n,N)& \leq \sum_{j=n+1}^N \mean{ \ti{b}_{1,\omega}(j)^2 } \Vert \calU(j)\Vert ^2 \mean{ \Vert \bfd^{+,k}_\omega(j)\Vert ^2} 
 \\
&\leq C \sum_{j=1}^\infty \mean{\ti{b}_{1,\omega}(j)^2}  \Vert \calU(j)\Vert ^2 =D^I < \infty.
\end{split}
\eeq
Similarly for $R^{II}(n,N)$, \eqref{decay-condition-b2} says
\beq \label{bound-on-RII}
\begin{split}
R^{II}(n,N) &\leq \sum_{j=n+1}^N \mean{ \ti{b}_{2,\omega}(j)^2 } \Vert \calV(j)\Vert ^2 \mean{\Vert \bfd^{+,k}_\omega(j)\Vert ^2} 
 \\
&\leq C \sum_{j=1}^\infty \mean{ \ti{b}_{2,\omega}(j)^2 } \Vert \calV(j)\Vert ^2 =D^{II}< \infty.
\end{split}
\eeq
and for $R^{III}(n,N)$, \eqref{decay-condition-b2b3} implies
\beq \label{bound-on-RIII}
\begin{split}
R^{III}(n,N) &\leq \sum_{j=n+1}^N \mean{ |\ti{b}_{2,\omega}(j)| \ti{b}_{3,\omega}(j) }  \Vert \calV(j) \Vert_{\rm HS} 
\Vert \calW(j) \Vert_{\rm HS} \mean{\Vert \bfd^{+,k}_\omega(j) \Vert^2 }  \\
&\leq C \sum_{j=1}^\infty \mean{ |\ti{b}_{2,\omega}(j)| \ti{b}_{3,\omega}(j) }  \Vert \calV(j) \Vert_{\rm HS} 
\Vert \calW(j) \Vert_{\rm HS} =D^{III}< \infty.
\end{split}
\eeq
The procedure we apply to $R^{IV}$ is a little more involved. Applying H\"older's inequality (in the
second inequality below) and then using independence we get 
\begin{align} \label{separate-RIV}
& R^{IV}(n,N) \notag \\
&=\mean{\left( \sum_{j=n+1}^N \left( \ti{b}_{3,\omega}(j)
\left(w_{11}(j)d^{+,k}_{1,\omega}(j)+w_{12}(j)d^{+,k}_{2,\omega}(j)\right) \right) \right)^2 } \notag \\ 
& \quad +\mean{\left( \sum_{j=n+1}^N \left( \ti{b}_{3,\omega}(j)
\left(w_{21}(j)d^{+,k}_{1,\omega}(j)+w_{22}(j)d^{+,k}_{2,\omega}(j)\right) \right) \right)^2 } \notag  \\
&\leq 2\sum_{i,j=n+1}^N \mean{ \ti{b}_{3,\omega}(j) \ti{b}_{3,\omega}(i) \Vert \bfd^{+,k}_\omega(j) \Vert 
\Vert \bfd^{+,k}_\omega(i) \Vert } \notag \\
&\quad \times \Vert \calW(j) \Vert_{\rm HS} \Vert \calW(i) \Vert_{\rm HS} \notag \\
& \leq 2\sum_{i,j=n+1}^N \mean{ \ti{b}_{3,\omega}(j)^2 \Vert \bfd^{+,k}_\omega(j) \Vert^2 } ^{1/2}
\mean{\ti{b}_{3,\omega}(i)^2 \Vert \bfd^{+,k}_\omega(i) \Vert^2 }^{1/2} \notag \\
& \quad \times \Vert \calW(j) \Vert_{\rm HS} \Vert \calW(i) \Vert_{\rm HS} \notag \\
&=2\sum_{i,j=n+1}^N \mean{\ti{b}_{3,\omega}(j)^2}^{1/2}
\mean{\Vert \bfd^{+,k}_\omega(j) \Vert^2 } ^{1/2}
\mean{\ti{b}_{3,\omega}(i)^2 }^{1/2}
\mean{ \Vert \bfd^{+,k}_\omega(i) \Vert^2 }^{1/2} \notag \\
&\quad \times \Vert \calW(j) \Vert_{\rm HS} \Vert \calW(i) \Vert_{\rm HS}  \notag \\
&=2 \Big( \sum_{j=n+1}^N \mean{ \ti{b}_{3,\omega}(j)^2 }^{1/2}
\mean{\Vert \bfd^{+,k}_\omega(j) \Vert^2 } ^{1/2}\Vert \calW(j) \Vert_{\rm HS} \Big)^2 \notag \\
&\leq 2C \Big( \sum_{j=n+1}^N \mean{ \ti{b}_{3,\omega}(j)^2 }^{1/2}\Vert \calW(j) \Vert_{\rm HS} \Big)^2 \notag \\
&\leq 2C \Big( \sum_{j=1}^\infty \mean{ \ti{b}_{3,\omega}(j)^2 }^{1/2}\Vert \calW(j) \Vert_{\rm HS} \Big)^2 =D^{IV} < \infty
\end{align}
by \eqref{decay-condition-b3}.
We see, therefore, that 
\beq \no
\begin{split}
\mean{ \Vert \bfd^{+,(k+1)}_\omega(n) \Vert^2 }&\leq
\liminf_{N \rightarrow \infty} \big( R^I(n,N)+R^{II}(n,N) \\ &+R^{III}(n,N)+R^{IV}(n,N) \big) 
\end{split}
\eeq
is bounded in $n$. It follows that $\bfd^{+,k}_\omega(n)$ is defined for all $k,n$ and $P$-almost every $\omega$.

The estimates above also imply that there exists a constant $C_0$ such that
\beq \label{recursion2}
\begin{split}
\mean{ \Vert \bfd^{+,k}_\omega(n) \Vert^2 } & 
\leq C_0 \big( \sup_{j>n} \mean{ \Vert \bfd^{+,(k-1)}_\omega(j) \Vert^2 } \big) \\ 
& \quad \times \Big( \sum_{j=n+1}^\infty \mean{\ti{b}_{1,\omega}(j)^2 } \Vert \calU(j)\Vert ^2+
\sum_{j=n+1}^\infty \mean{\ti{b}_{2,\omega}(j)^2 } \Vert \calV(j)\Vert ^2 \\
& \quad +\sum_{j=n+1}^\infty \mean{|\ti{b}_{2,\omega}(j)| \ti{b}_{3,\omega}(j) }  \Vert \calV(j) \Vert_{\rm HS} 
\Vert \calW(j) \Vert_{\rm HS} \\ 
& \quad +\big( \sum_{j=n+1}^\infty \mean{ \ti{b}_{3,\omega}(j)^2 }^{1/2}\Vert \calW(j) \Vert_{\rm HS} \big)^2 \Big)
\end{split}
\eeq
which, as in the proof of the previous lemma, implies in turn that for large enough $n$, 
\beq \label{smallness-of-mean-for-a+k}
\mean{ \Vert \bfd^{+,k}_\omega(n) \Vert }^2 \leq \mean{ \Vert \bfd^{+,k}_\omega(n) \Vert^2 } \leq \left( \frac{1}{4} \right)^k 
\eeq
and therefore that for such $n$,
\beq \no
\bfd^+_\omega(n)=\sum_{k=0}^\infty \bfd^{+,k}_\omega(n)
\eeq
converges almost surely.

The next step is to show
\begin{align} \label{summation-order2}
&\sum_{k=0}^\infty \sum_{j=n+1}^\infty \left( \ti{b}_{1,\omega}(j)\calU(j)+\ti{b}_{2,\omega}(j)\calV(j)+\ti{b}_{3,\omega}(j)\calW(j) \right) 
\bfd^{+,k}_\omega(j) \notag \\
&=\sum_{j=n+1}^\infty \sum_{k=0}^\infty \left( \ti{b}_{1,\omega}(j)\calU(j)+\ti{b}_{2,\omega}(j)\calV(j)+\ti{b}_{3,\omega}(j)\calW(j) \right) 
\bfd^{+,k}_\omega(j) 
\end{align}
 from which the first half of the theorem (\eqref{asymptotic_amplitude12} and \eqref{asymptotic_amplitude22}), will follow.
Note, first, that from \eqref{smallness-of-mean-for-a+k} and \eqref{decay-condition-b3}, $P$-a.s.\
\beq \label{summation-order-partial1}
\sum_{k=0}^\infty \sum_{j=n+1}^\infty \left\| \ti{b}_{3,\omega}(j)\calW(j)\bfd^{+,k}_\omega(j) \right\| <\infty
\eeq
so that we have
\beq \label{summation-order-partial2}
\sum_{k=0}^\infty \sum_{j=n+1}^\infty \ti{b}_{3,\omega}(j)\calW(j)\bfd^{+,k}_\omega(j)=
\sum_{j=n+1}^\infty \sum_{k=0}^\infty \ti{b}_{3,\omega}(j)\calW(j)\bfd^{+,k}_\omega(j)
\eeq
with probability one. Thus we are only left with showing
\begin{align} \label{summation-order-partial3}
&\sum_{k=0}^\infty \sum_{j=n+1}^\infty \left( \ti{b}_{1,\omega}(j)\calU(j)+\ti{b}_{2,\omega}(j)\calV(j) \right) 
\bfd^{+,k}_\omega(j) \notag \\
&=\sum_{j=n+1}^\infty \sum_{k=0}^\infty \left( \ti{b}_{1,\omega}(j)\calU(j)+\ti{b}_{2,\omega}(j)\calV(j) \right) 
\bfd^{+,k}_\omega(j). 
\end{align}
The proof of \eqref{summation-order-partial3} is precisely the same as the corresponding step in the proof of 
Lemma~\ref{random_variation}. It is therefore omitted from the argumentation.

Obviously, to show \eqref{asymptotic_amplitude32} and \eqref{asymptotic_amplitude42} one simply follows the exact same 
procedure outlined above, with a different initial vector, so we are done. 
\end{proof}

\section{Proof of Theorems \ref{ac_stability} and 
\ref{singular_stability}}

In this section we present the proofs of Theorems \ref{ac_stability} 
and \ref{singular_stability}.
\begin{proof}[Proof of Theorem~\ref{ac_stability}]
There are two cases to consider:
\begin{itemize}
\item{case 1.} Assume that there exists a subsequence $a(n_j) \rightarrow 0$ as $j \rightarrow \infty$. Then we may choose
a summable subsequence $a(k_j)$. Define 
\begin{displaymath} \no
a^1(k)= \left\{ \begin{array}{ll}
a(k_j) & \textrm{if } k=k_j \textrm{ for some } j \\
0 & \textrm{otherwise}
\end{array} \right.
\end{displaymath}
and let $\ti{J}_0=J_0-J \left( \se{a^1},\szer \right)$. Then $\ti{J}_0$ is a direct sum of finite rank operators so 
that its spectrum is pure point. Since $J_0$ is a trace class perturbation of $\ti{J}_0$, it follows that $\ti{J}_0$
has no absolutely continuous spectrum. If $\Gamma$ is empty as well, then we are done. Otherwise, \eqref{Gamma-def} 
implies that $\ti{a}_\omega(n) \rightarrow 0$ as $n \rightarrow \infty$ almost surely. 
This is because $ \Vert T^E(n) \Vert \geq C \min(1,a(n)^{-(1/2)})$ for some universal constant $C$.
Therefore, repeating the argument
above for $J_\omega$, we find that it has no absolutely continuous spectrum as well. Thus \eqref{ac_conclusion} 
follows for this case.

\item{case 2.} There is a constant $c_0>0$ such that $a(n)>c_0$ for any $n$. Then it follows from 
\eqref{condition_on_alpha}, that 
\beq \label{boundedness-of-offdiagonals}
\frac{1}{a(n)+\ti{a}_\omega(n)} < (\delta c_0)^{-1}
\eeq
for all $n$ and $\omega$.

Fix $E\in\Gamma$ and let $\{N_j\}_{j=1}^\infty$ be a sequence for which
\beq \no
\lim_{j \rightarrow \infty} \left( \frac{1}{N_j}\sum_{n=1}^{N_j}\frac{1}{a(n)} \right)>0.
\eeq
We want to show that with probability 
one, \begin{eqnarray} \label{local-condition}
& \liminf_{j\rightarrow\infty}
\frac{1}{N_j}\sum_{n=1}^{N_j} \parallel T^E_0(n)\parallel ^2<\infty \nonumber \\
& \Longleftrightarrow \nonumber \\
& \liminf_{j\rightarrow\infty}
\frac{1}{N_j}\sum_{n=1}^{N_j} \parallel T^E_\omega(n)\parallel ^2<\infty\
\end{eqnarray}
(where we use the notation introduced in the beginning of the previous 
section). Then, by Fubini, it will follow that there exists a set of full  
$P$-measure of realizations of the perturbation, such that for 
Lebesgue-a.e.\ energy in $\Gamma$, the asymptotic properties of the 
transfer matrices (in the above sense) remain the same. By 
Proposition~\ref{last-simon}, this implies \eqref{ac_conclusion} (since $\ti{a}_\omega(n) \rightarrow 0$ almost surely).

The one-step transfer matrices have the form \beq \no
S^E_\omega(n) = \left( \begin{array}{cc}
\frac{E-b(n)-\ti{b}_\omega(n)}{a(n)+\ti{a}_\omega(n)} & 
\frac{-a(n-1)-\ti{a}_\omega(n-1)}{a(n)+\ti{a}_\omega(n)} \\
1                       &  0
\end{array} \right)
\eeq
where $a(0)+\ti{a}_\omega(0) \equiv 1$. Due to the 
$\frac{a(n-1)+\ti{a}_\omega(n-1)}{a(n)+\ti{a}_\omega(n)}$ term, these 
matrices are not independent. A crucial ingredient in the proof of Lemma~\ref{random_variation2} 
is the independence of the matrices.
We begin, therefore, with a modification to 
these matrices following \cite{kala}:
Define
\beq \label{Komega}
K_\omega(n)=\left( \begin{array}{cc}
1 & 0 \\
0 & a(n)+\ti{a}_\omega(n) 
\end{array} \right).
\eeq
Then 
\begin{align} \label{Stilde}
&K_\omega(n)S^E_\omega(n)K_\omega(n-1)^{-1}=\left( \begin{array}{cc}
\frac{E-b(n)-\ti{b}_\omega(n)}{a(n)+\ti{a}_\omega(n)} & 
\frac{-1}{a(n)+\ti{a}_\omega(n)} \\
a(n)+\ti{a}_\omega(n) & 0
\end{array} \right) \notag \\
&\equiv
\ti{S}^E_\omega(n).
\end{align}
Note that $\ti{S}^E_\omega(n)$ are independent and unimodular. One may 
now define
\beq \label{Ttilde}
\begin{split}
\ti{T}^E_\omega(n) \equiv 
\ti{S}^E_\omega(n)\cdot\ldots \cdot \ti{S}^E_\omega(1)&=
K_\omega(n)T^E_\omega(n)K^E_\omega(0)^{-1} \\
&=K_\omega(n)T^E_\omega(n).
\end{split}
\eeq
Define $\ti{\calD}^E_\omega(n)$ through
\beq \label{Atilde}
\ti{T}^E_\omega(n)=\ti{T}^E_0(n)\ti{\calD}^E_\omega(n).
\eeq
Then
\beq \no
T^E_\omega(n)=K_\omega(n)^{-1}K_0(n)T^E_0(n)\ti{\calD}^E_\omega(n),
\eeq
so, using \eqref{condition_on_alpha}, in order to show \eqref{local-condition} 
almost surely, it suffices to show that $\ti{\calD}^E_\omega(n)$ 
converge to a limit with probability one. 

 From \eqref{Atilde} it follows that
\beq \label{Atilde-recursion}
\ti{\calD}^E_\omega(n-1)=(I+\ti{\calU}^E_\omega(n))\ti{\calD}^E_\omega(n),
\eeq
where
\begin{align} \label{Utilde-definition}
&\ti{\calU}^E_\omega(n) \notag \\
&=\ti{T}^E_0(n)^{-1}(\ti{S}^E_0(n)\ti{S}^E_\omega(n)^{-1}-I) \ti{T}^E_0(n) \notag \\
&=\ti{T}^E_0(n)^{-1} \left( \begin{array}{cc}
\frac{\ti{a}_\omega(n)}{a(n)} & \frac{\ti{b}_\omega(n)}{a(n)(a(n)+\ti{a}_\omega(n))} \\
0 & -\frac{\ti{a}_\omega(n)}{a(n)+\ti{a}_\omega(n)}
\end{array} \right)\ti{T}^E_0(n) \notag \\
&=\frac{\ti{a}_\omega(n)}{a(n)}\ti{T}^E_0(n)^{-1} \left( \begin{array}{cc}
1 & 0 \\
0 & 0
\end{array} \right)\ti{T}^E_0(n) \notag \\
& \quad +\frac{\ti{b}_\omega(n)}{a(n)(a(n)+\ti{a}_\omega(n))}  \ti{T}^E_0(n)^{-1} \left( \begin{array}{cc}
0 & 1 \\
0 & 0
\end{array} \right)\ti{T}^E_0(n) \notag \\
& \quad +\frac{\ti{a}_\omega(n)}{a(n)+\ti{a}_\omega(n)} \ti{T}^E_0(n)^{-1} \left( \begin{array}{cc}
0 & 0 \\
0 & -1
\end{array} \right)\ti{T}^E_0(n). 
\end{align}
Writing 
\beq \no
\frac{\ti{a}_\omega(n)}{a(n)+\ti{a}_\omega(n)}=\frac{\ti{a}_\omega(n)}{a(n)}
-\frac{(\ti{a}_\omega(n))^2}{a(n)(a(n)+\ti{a}_\omega(n))},
\eeq 
and denoting 
\beq \label{U-definition}
\calV^E(n)=\ti{T}^E_0(n)^{-1} \left( \begin{array}{cc}
1 & 0 \\
0 & -1
\end{array} \right)\ti{T}^E_0(n),
\eeq
\beq \label{V-definition}
\calU^E(n)=\ti{T}^E_0(n)^{-1} \left( \begin{array}{cc}
0 & 1 \\
0 & 0
\end{array} \right)\ti{T}^E_0(n),
\eeq
\beq \label{W-definition}
\calW^E(n)=\ti{T}^E_0(n)^{-1} \left( \begin{array}{cc}
0 & 0 \\
0 & 1
\end{array} \right)\ti{T}^E_0(n),
\eeq
we get
\beq \label{Utilde-definition2}
\begin{split}
\ti{\calU}^E_\omega(n)&=\frac{\ti{a}_\omega(n)}{a(n)}\calV^E(n) 
+\frac{\ti{b}_\omega(n)}{a(n)(a(n)+\ti{a}_\omega(n))} \calU^E(n)+ \\
& \quad +\frac{(\ti{a}_\omega(n))^2}{a(n)(a(n)+\ti{a}_\omega(n))} \calW^E(n).
\end{split}
\eeq
Since the $\ti{T}^E_0(n)$ are unimodular (so they have norm equal to their inverses') it follows that 
\beq \label{bound-on-U}
\Vert \calU^E(n)\Vert \leq \Vert \ti{T}^E_0(n)\Vert ^2
\eeq 
and the same holds for $\calV^E(n)$ and $\calW^E(n)$.
Furthermore, \eqref{Gamma-def} implies that
\begin{align} \no
&\sum_{n=1}^\infty 
\left(\mean{ \ti{a}_\omega(n)^4 } ^{1/2}+\mean{ \ti{b}_\omega(n)^2 } \right) \Vert \ti{T}^E_0(n)\Vert ^4 
\notag \\
&\leq \sum_{n=1}^\infty 
\left(\mean{ \ti{a}_\omega(n)^4 } ^{1/2}+\mean{ \ti{b}_\omega(n)^2 } \right) \Vert K_0(n) \Vert^4
\Vert T^E_0(n)\Vert ^4 \notag \\
&\leq C\sum_{n=1}^\infty 
\left(\mean{ \ti{a}_\omega(n)^4 } ^{1/2}+\mean{ \ti{b}_\omega(n)^2 } \right)(a(n)+1)^4 
\Vert T^E_0(n)\Vert ^4 
<\infty.
\end{align}
(Here $C$ is some universal constant.)
Thus, from \eqref{boundedness-of-offdiagonals} and from \eqref{Gamma-def}, it follows that
\begin{align} \label{condition4-1}
&\sum_{n=1}^\infty \mean{ \frac{\ti{b}_\omega(n)^2}{a(n)^2(a(n)+\ti{a}_\omega(n))^2}} \Vert \calU^E(n) \Vert_{\rm HS}^2 
\notag \\
&\leq C_1 \sum_{n=1}^\infty \mean{ \ti{b}_\omega(n)^2 } \Vert \ti{T}^E_0(n) \Vert^4 <\infty
\end{align} 
and
\begin{align} \label{condition4-2}
&\sum_{n=1}^\infty \mean{ \frac{\ti{a}_\omega(n)^2}{a(n)^2}} \Vert \calV^E(n) \Vert_{\rm HS}^2 
\notag \\
&\leq C_2 \sum_{n=1}^\infty \mean{\ti{a}_\omega(n)^4 }^{1/2} \Vert \ti{T}^E_0(n) \Vert^4 <\infty
\end{align}
for some constants $C_1,\ C_2>0$. Also, since $\Vert \ti{T}^E_0(n) \Vert \geq 1$ (recall these are unimodular),
\begin{align} \label{condition5}
&\sum_{n=1}^\infty \mean{\frac{|\ti{a}_\omega(n)|\ti{a}_\omega(n)^2}{a(n)^2(a(n)+\ti{a}_\omega(n))}} 
\Vert \calV^E(n) \Vert_{\rm HS} \Vert \calW^E(n) \Vert_{\rm HS} \notag \\
& \leq C_3 \sum_{n=1}^\infty \mean{ |\ti{a}_\omega(n)|^3 } 
\Vert \ti{T}^E_0(n) \Vert^4  \notag \\
& \leq C_3 \sum_{n=1}^\infty \mean{ |\ti{a}_\omega(n)|^4 }^{3/4} 
\Vert \ti{T}^E_0(n) \Vert^4  \notag \\
&=C_3 \sum_{n=1}^\infty \big( \mean{ |\ti{a}_\omega(n)|^4 }^{1/2} 
\Vert \ti{T}^E_0(n) \Vert^4 \big) ^{3/2} < \infty
\end{align}
where $C_3>0$ is a constant.
Finally, again using the fact that $\Vert \ti{T}^E_0(n) \Vert \geq 1$,
we get, for some constant $C_4>0$,
\begin{align} \label{condition6}
&\sum_{n=1}^\infty \mean{ \frac{\ti{a}_\omega(n)^4}{a(n)^2(a(n)+\ti{a}_\omega(n))^2} }^{1/2} \Vert \calW^E(n) \Vert_{\rm HS} 
\notag \\
& \leq C_4 \sum_{n=1}^\infty \mean{\ti{a}_\omega(n)^4 }^{1/2} \Vert \ti{T}^E_0(n) \Vert^2 \notag \\
& \leq C_4 \sum_{n=1}^\infty \mean{\ti{a}_\omega(n)^4 }^{1/2} \Vert \ti{T}^E_0(n) \Vert^4<\infty.
\end{align}
Thus, we see that the conditions of Lemma~\ref{random_variation2} are satisfied, with 
$\ti{b}_1 \equiv \frac{\ti{b}}{a(a+\ti{a})}$, $\ti{b}_2 \equiv \frac{\ti{a}}{a}$, 
$\ti{b}_3 \equiv \frac{\ti{a}^2}{a(a+\ti{a})}$ and the obvious correspondence for the matrices. This implies that with
probability one, the matrices $\ti{\calD}^E_\omega(n)$ converge to the identity matrix. 
As explained above, this finishes the proof of case 2 and therefore completes the proof of the theorem.
\end{itemize} 
\end{proof}
\begin{remark}
It is important to note that for a perturbation along the diagonal alone (that is - for the case of 
$\ti{a}_\omega(n) \equiv 0$), a much shorter proof can be provided:
Note that in this case, $\ti{\calU}^E_\omega(n)$ of \eqref{Atilde-recursion} reduces to 
\beq \label{Utilde-reduction}
\ti{\calU}^E(\omega,n)=\frac{\ti{b}_\omega(n)}{a(n)^2}\calU^E(n)
\eeq
so that
\beq \label{zero-mean-in-recursion}
\mean{ \ti{\calU}^E_\omega(n)}=\frac{\mean{ \ti{b}_\omega(n) }}{a(n)^2}\calU^E(n) =0.
\eeq 
It is clear also that $\det \left( I+\ti{\calU}^E_\omega(n) \right)=1$ so that we have also
\beq \label{martingale-property-for-A}
\mean{ \left( I+\ti{U}^E_\omega,n)\right)^{-1} } = I.
\eeq
Thus, from equation \eqref{Atilde-recursion}, and since 
$\ti{\calU}^E_\omega(n)$ is a function of the perturbing potential at the point $n$ alone,
it follows that $\ti{\calD}^E_\omega(n)$ is a matrix-valued martingale. 
$\Gamma$ is precisely the set where this martingale is bounded, so the theorem follows from the martingale convergence 
theorem. The extra work we do (in Lemma~\ref{random_variation2}) is due to the term 
$\frac{\ti{a}_\omega(n)^2}{a(n)(a(n)+\ti{a}_\omega(n))}$ which does not have zero mean but is $\ell^1$ almost surely.
\end{remark}
\begin{proof}[Proof of Theorem~\ref{singular_stability}]
Recall the definition of $\vphi_{1,\theta}$ and $\vphi_{2,\theta}$ 
(\eqref{phi1theta} and \eqref{phi2theta}) for a given potential 
$\se{b}$. Note that for the operator $H_\theta$ the transfer matrices 
have the form 
\beq \label{transfer-solutions}
T^E(n)=\left( \begin{array}{cc}
\vphi_{1,\theta}(n+1) & \vphi_{2,\theta}(n+1) \\
\vphi_{1,\theta}(n)   & \vphi_{2,\theta}(n)  
\end{array} \right).
\eeq

For the perturbing random potential $\se{\ti{b}_\omega}$, and given 
the characterization \eqref{singularity_condition} of $\Sigma_{\alphas}$, 
it is obvious that in
order to prove stability of the local Hausdorff dimension for a given
energy - $E$, it suffices to show that a.s.\ there exist two solutions
$\psi_{1,\omega}(n)$ and $\psi_{2,\omega}(n)$ of 
\beq \label{perturbed-difference}
\psi(n+1)+\psi(n-1)+(b(n)+\ti{b}_\omega(n))\psi(n)=E\psi(n)
\eeq
that satisfy 
\begin{eqnarray} \label{asymptotic_equality1}
\lim_{L \rightarrow \infty} 
\frac {\parallel \psi_{1,\omega} \parallel _L}{\parallel \vphi_{1,\theta(E)}\parallel _L}
& = & 1 \\
\label{asymptotic_equality2}
\lim_{L \rightarrow \infty} 
\frac {\parallel \psi_{2,\omega} \parallel _L}{\parallel \vphi_{2,\theta(E)}\parallel _L}
& = & 1.
\end{eqnarray}
We shall prove that relations \eqref{asymptotic_equality1} and 
\eqref{asymptotic_equality2} hold almost surely, for every energy in the 
set $\Lambda$. A simple application of Fubini's theorem (as in the 
previous proof), then yields the inclusion in the theorem up to a set of 
measure zero, for any fixed measure.

Fix $E \in \Lambda$. It isn't difficult to see (see Lemma~4.3  in 
\cite{kls}) that for any $\veps>0$ there exist $\veps$-dependent constants $C_1,\ C_2,\ 
C_3,\ 
C_4,$ so that for large $N$
\beq \label{bound_on_phi1}
C_2N^{1-\frac{1}{2\beta(E)}-\veps} \leq \parallel \vphi_{1,\theta(E)}\parallel _N \leq C_1 
N^{1/2+\veps}
\eeq
\beq \label{bound_on_phi2}
C_4N^{1/2-\veps} \leq \parallel \vphi_{2,\theta(E)}\parallel _N \leq 
C_3N^{\frac{1}{2\beta(E)}+\veps}.
\eeq
(Recall that $\beta(E)>0$.)
Now, applying Lemma~\ref{random_variation}, with 
\beq \label{singularU}
\begin{split}
\calU(n)&=\left( 
\begin{array}{cc}
\vphi_{1,\theta(E)}(n+1) & \vphi_{2,\theta(E)}(n+1) \\
\vphi_{1,\theta(E)}(n)   & \vphi_{2,\theta(E)}(n)
\end{array} \right)^{-1}\cdot \\
& \quad \cdot \left(
\begin{array}{cc}
0 & 1 \\
0  & 0
\end{array} \right)
\left(
\begin{array}{cc}
\vphi_{1,\theta(E)}(n+1) & \vphi_{2,\theta(E)}(n+1) \\
\vphi_{1,\theta(E)}(n)   & \vphi_{2,\theta(E)}(n)
\end{array} \right),
\end{split}
\eeq
$\ti{b}_\omega(n),$ and $f_+(n)=n^{\ti{\eta}}$ (recall \eqref{Lambda}), we see that with probability 
one, there exist sequences $\bfd^+_\omega(n) \equiv \left ( \begin{array} {c} d^+_{1,\omega}(n) \\
d^+_{2,\omega}(n)
\end{array} \right)$ 
and 
$\bfd^-_\omega(n) \equiv \left ( \begin{array} {c} d^-_{1,\omega}(n) \\
d^-_{2,\omega}(n)
\end{array} \right)$, 
that solve \eqref{conclusion_for_a} and satisfy 
\eqref{asymptotic_amplitude1}--\eqref{asymptotic_amplitude4}.
Let 
\beq \label{psi1}
\begin{split}
\psi_{1,\omega}(n) &= d_{1,\omega}^-(n)\vphi_{1,\theta(E)}(n)+
d_{2,\omega}^-(n)\vphi_{2,\theta(E)}(n) \\
&=d_{1,\omega}^-(n-1)\vphi_{1,\theta(E)}(n)+
d_{2,\omega}^-(n-1)\vphi_{2,\theta(E)}(n)
\end{split}
\eeq
\beq \label{psi2}
\begin{split}
\psi_{2,\omega}(n) &= d_{1,\omega}^+(n)\vphi_{1,\theta(E)}(n)+
d_{2,\omega}^+(n)\vphi_{2,\theta(E)}(n) \\
&= d_{1,\omega}^+(n-1)\vphi_{1,\theta(E)}(n)+
d_{2,\omega}^+(n-1)\vphi_{2,\theta(E)}(n).
\end{split}
\eeq
(The last equality in each equation follows from \eqref{conclusion_for_a} with the above definition for $\calU(n)$.)
Using 
\beq \no
\vphi_{1,\theta(E)}(n)\vphi_{2,\theta(E)}(n-1)-\vphi_{1,\theta(E)}(n-1)\vphi_{2,\theta(E)}(n)=\det T^E(n-1)=1
\eeq 
and \eqref{conclusion_for_a}, we get that
\begin{align} \no
&d_{1,\omega}^\pm(n-1)\vphi_{1,\theta(E)}(n-1)+d_{2,\omega}^\pm(n-1)\vphi_{2,\theta(E)}(n-1) 
\notag \\
&=d_{1,\omega}^\pm(n)\vphi_{1,\theta(E)}(n-1)+d_{2,\omega}^\pm(n)\vphi_{2,\theta(E)}(n-1) \notag \\
& \quad -\ti{b}_\omega(n)(d_{1,\omega}^\pm(n)\vphi_{1,\theta(E)}(n)+d_{2,\omega}^\pm(n)\vphi_{2,\theta(E)}(n)). \notag
\end{align}
Thus we see that
\begin{align} \no
&d_{1,\omega}^\pm(n) \vphi_{1,\theta(E)}(n+1)+d_{2,\omega}^\pm(n) \vphi_{2,\theta(E)}(n+1) \notag \\
& \quad +d_{1,\omega}^\pm(n-1) \vphi_{1,\theta(E)}(n-1)+d_{2,\omega}^\pm(n-1) \vphi_{1,\theta(E)}(n-1) \notag \\
& \quad +\big(b(n)+\ti{b}_\omega(n) \big)
\big( d_{1,\omega}^\pm(n) \vphi_{1,\theta(E)}(n)+d_{2,\omega}^\pm(n) \vphi_{2,\theta(E)}(n) \big) \notag \\
&=E(d_{1,\omega}^\pm(n) \vphi_{1,\theta(E)}(n)+d_{2,\omega}^\pm(n) \vphi_{2,\theta(E)}(n)) \notag
\end{align}
which means that $\se{\psi_{1,\omega}}$ and $\se{\psi_{2,\omega}}$ solve \eqref{perturbed-difference}.
 
Now, from 
\beq \no
\begin{split}
\left\vert \frac{\Vert \psi_{2,\omega} \Vert_L}{\Vert \vphi_{2,\theta(E)}\Vert_L}-1 \right\vert &\leq
\frac{\Vert \psi_{2,\omega}-\vphi_{2,\theta(E)} \Vert_L}{\Vert \vphi_{2,\theta(E)}\Vert_L} \\
&=\frac{\Vert d_{1,\omega}^+\vphi_{1,\theta(E)}+(d_{2,\omega}^+-1)\vphi_{2,\theta(E)} \Vert_L}{\Vert 
\vphi_{2,\theta(E)}\Vert_L}
 \\
&\leq \sqrt{2} \left( \frac{\Vert d_{1,\omega}^+\vphi_{1,\theta(E)} \Vert_L}{\Vert \vphi_{2,\theta(E)}\Vert_L}+
\frac{\Vert (d_{2,\omega}^+-1)\vphi_{2,\theta(E)} \Vert_L}{\Vert \vphi_{2,\theta(E)}\Vert_L} \right) 
\end{split}
\eeq
it follows immediately that \eqref{asymptotic_equality2} holds, (since $\vphi_{1,\theta(E)}$ is subordinate). 
For \eqref{asymptotic_equality1} write, similarly,
\beq \no
\begin{split}
\left\vert \frac{\Vert \psi_{1,\omega} \Vert_L}{\Vert \vphi_{1,\theta(E)}\Vert_L}-1 \right\vert & \leq
\frac{\Vert \psi_{1,\omega}-\vphi_{1,\theta(E)} \Vert_L}{\Vert \vphi_{1,\theta(E)}\Vert_L} \\
&=\frac{\Vert (d_{1,\omega}^--1)\vphi_{1,\theta(E)}+d_{2,\omega}^-\vphi_{2,\theta(E)} \Vert_L}{\Vert 
\vphi_{1,\theta(E)}\Vert_L}
\\
&\leq \sqrt{2} \left( \frac{\Vert (d_{1,\omega}^--1)\vphi_{1,\theta(E)} \Vert_L}{\Vert \vphi_{1,\theta(E)}\Vert_L}+
\frac{\Vert d_{2,\omega}^-\vphi_{2,\theta(E)} \Vert_L}{\Vert \vphi_{1,\theta(E)}\Vert_L} \right). 
\end{split}
\eeq
That the first term on the left hand side converges to zero is immediate (recall that $\beta(E)>0$ so that 
$\lim_{L \rightarrow \infty} \Vert \vphi_{1,\theta(E)} \Vert_L=\infty$). For the second term we note, first, that for 
any $\veps>0$, there is an $L_0$ such that for any $n \geq L_0$ 
$|d_{2,\omega}^-(n)|<\frac{\veps}{n^{\ti{\eta}}}$. Second, from \eqref{bound_on_phi1} and \eqref{bound_on_phi2}, we get that
there exists a constant $D>0$ for which 
\beq \no
\frac{\Vert \vphi_{2,\theta(E)} \Vert_L}{\Vert \vphi_{1,\theta(E)}\Vert_L} \leq DL^{\eta(E)}.
\eeq
Now, using summation by parts,
\begin{align} \no
&\frac{\left\| d_{2,\omega}^-\vphi_{2,\theta(E)} \right\|_L}{\left\| \vphi_{1,\theta(E)}\right\|_L} \notag \\
&=\frac{\left( \sum_{n=1}^L d_{2,\omega}^-(n)^2|\vphi_{2,\theta(E)}(n)|^2 \right)^{1/2}}{\left\| \vphi_{1,\theta(E)}\right\|_L} \notag 
\\
&\leq \frac{\left\| d_{2,\omega}^-\vphi_{2,\theta(E)} \right\|_{L_0}}{\left\| \vphi_{1,\theta(E)}\right\|_L}+
\frac{\left( \sum_{n=1}^L \frac{\veps^2}{n^{2\ti{\eta}}}|\vphi_{2,\theta(E)}(n)|^2 \right)^{1/2}}{\left\| \vphi_{1,\theta(E)}\right\|_L} \notag \\  
&\leq \frac{\left\| d_{2,\omega}^-\vphi_{2,\theta(E)} \right\|_{L_0}}{\left\| \vphi_{1,\theta(E)}\right\|_L}+
\frac{\veps}{L^{\ti{\eta}}} \frac{\left\| \vphi_{2,\theta(E)} \right\|_L}{\left\| \vphi_{1,\theta(E)} \right\|_L}\notag \\
& \quad +\frac{\left( \sum_{n=1}^L \left( (\frac{\veps^2}{n^{2\ti{\eta}}}-\frac{\veps^2}{(n+1)^{2\ti{\eta}}}) 
\sum_{j=1}^n |\vphi_{2,\theta(E)}(j)|^2 \right) \right)^{1/2}}{\left\| \vphi_{1,\theta(E)} \right\|_L} \notag \\
&\leq \frac{\left\| d_{2,\omega}^-\vphi_{2,\theta(E)} \right\|_{L_0}}{\left\| \vphi_{1,\theta(E)}\right\|_L}+
D \cdot \frac{\veps}{L^{\ti{\eta}}} L^{\eta(E)} \notag \\
& \quad +D \cdot \veps \frac{\left( \sum_{n=1}^L\frac{1}{n^{2\ti{\eta}+1}}n^{2\eta(E)} 
\sum_{j=1}^n |\vphi_{1,\theta(E)}(j)|^2 \right)^{1/2}}{\left\| \vphi_{1,\theta(E)} \right\|_L} \notag \\
&\leq \frac{\left\| d_{2,\omega}^-\vphi_{2,\theta(E)} \right\|_{L_0}}{\left\| \vphi_{1,\theta(E)}\right\|_L}+
D \cdot \veps 
+D \cdot \veps \left( \sum_{n=1}^L\frac{1}{n^{1+2(\ti{\eta}-\eta(E))}}\right)^{1/2}, \notag 
\end{align}
so \eqref{asymptotic_equality1} follows from the fact that $\ti{\eta}>\eta(E)$. This finishes the proof of the theorem.
\end{proof}


\section{An Application of Theorem~\ref{singular_stability} to Sparse Potentials}

In this section, we present an application of Theorem~\ref{singular_stability} to one-dimensional Schr\"odin\-ger operators with 
sparse potentials, studied by Zlato\v s in \cite{zlatos}. The family 
$H_\theta(v,\gamma)$ of operators constructed there has a potential of the 
form
\beq \label{zlatos_potential}
b^{v,\gamma}(n)=\left\{ \begin{array}{cc}
v & n=n_j \equiv \gamma^j \textrm{ for some } j\geq1 \\
0 & \textrm{ otherwise.}
\end{array} \right.
\eeq
for some $v \neq 0$ and $\gamma>1$ an integer. 

We say that a measure $\mu$ 
has \emph{fractional Hausdorff dimension} in some interval $I$ if 
$\mu(I\cap\cdot)$ is $\alpha$-continuous and $(1-\alpha)$-singular for 
some $1>\alpha>0$. Considering potentials of the form 
\eqref{zlatos_potential}, Zlato\v s proves the following:
 
\begin{proposition}[Theorem~4.1 in \cite{zlatos}] \label{zlatos_thm}
For any closed interval of energies $I \subseteq (-2,2)$ there are $v_0>0$ 
and $\gamma_0 \in \bbN$ such that if $0<|v|<v_0$ and $\gamma \geq 
\gamma_0$ is an integer, then for Lebesgue-almost every $\theta$, the 
measure $\mu_\theta$, corresponding to the operator $H_\theta(v,\gamma)$, 
has fractional Hausdorff dimension in $I$. 
\end{proposition}
An important feature of operators with sparse potentials, is that the 
modulus of the solutions undergoes significant changes only near 
the 
points where the potential does not vanish (this is easily seen using EFGP 
transform - see \cite{efgp}). Thus, it is possible to obtain estimates on 
the pointwise behavior of the solutions looking at  points in the support of 
the potential. The proof of Proposition~\ref{zlatos_thm} goes 
through such estimates. 

It is actually shown there that, given $\gamma$ large enough and $v$ small enough, there exist constants 
$\beta_1<\beta_2<\frac{1}{2}$, depending only on $I,\ \gamma$ and $v$, and a set $I' \subseteq I$ of full Lebesgue measure,
such that for any $E \in I'$, equation \eqref{difference-schrodinger} with $b(n)=b^{v,\gamma}(n)$ has two solutions - 
$\vphi_1^E$ and $\vphi_2^E$ for which the following holds for sufficiently large $n$ and some small $\veps$:
\beq \label{zlatos-decay-for-phi1}
\big( |\vphi_1^E(n-1)|^2+|\vphi_1^E(n)|^2 \big)^{1/2} \leq n_j^{-\beta_1-\veps} \quad n_j<n \leq n_{j+1},
\eeq 
\beq \label{zlatos-decay-for-phi2}
n_j^{\beta_1} \leq \big( |\vphi_2^E(n-1)|^2+|\vphi_2^E(n)|^2 \big)^{1/2} \leq n_j^{\beta_2} \quad n_j<n \leq n_{j+1}.
\eeq 
 From this, using subordinacy theory \cite{jit-last} and the theory of rank one perturbations \cite{rankone}, Zlato\v s 
shows that for almost every boundary condition, $\theta$ the spectral measure $\mu_\theta$ is $(1-2\beta_2)$-continuous and
$(1-2\beta_1)$-singular on I. 

\eqref{zlatos-decay-for-phi1} and \eqref{zlatos-decay-for-phi2} provide us with a natural setting
to apply Theorem~\ref{singular_stability}. For the sake of simplicity and to make things explicit, we shall examine 
perturbing potentials of the form
\beq \label{perturbation-zlatos}
\ti{b}_{s,\omega}(n)=\frac{X_\omega(n)}{n^s}
\eeq
where $\se{X_\omega}$ are i.i.d.\ with a uniform distribution on an interval (say $[-1,1]$) and $s>0$ to
be specified later.

\begin{theorem}
Let $I\subseteq (-2,2)$ be a closed interval and assume that $H_\theta(v,\gamma)$ is an operator satisfying the 
requirements of Proposition~\ref{zlatos_thm} so that for a.e.\ $\theta$ its spectral measure is $\alpha$-continuous
and $(1-\alpha)$-singular on J, for some $0<\alpha<1$. Let 
\beq \label{beta1}
\beta_1=\beta_1(I,\gamma,v)
\eeq
and
\beq \label{beta2}
\beta_2=\beta_2(I,\gamma,v)
\eeq
be as in the discussion above, and let 
\beq \label{s-condition-zlatos}
s > \frac{4\beta_2}{1-2\beta_2}-2\beta_1+\frac{1}{2}.
\eeq
Then $P$-almost surely, the spectral measure of the random operator
\beq \label{perturbed-operator-zlaots}
H_\omega(v,\gamma,s)=H(v,\gamma)+\ti{b}_{s,\omega}
\eeq
is $\alpha$-continuous and $(1-\alpha)$-singular on I.
\end{theorem}
\begin{remark}
It is not hard to see that, using the result of Kiselev-Last-Simon \cite{kls} 
described in the introduction, one needs to demand 
\beq \label{kls_zlatos}
|\ti{b}(n)| \leq Cn^{-(s+\frac{1}{2})},
\eeq
in order to obtain this kind of stability.
\end{remark}

\begin{proof}
We want to apply the perturbation only to sites $2,3....$ at first, so denote by 
$\ti{b}^0_\omega(n)$ the sequence
\begin{displaymath} \no
\ti{b}^0_{s,\omega}(n)= \left\{ \begin{array}{ll}
\ti{b}_{s,\omega}(n) & \textrm{if } n>1 \\
0 & \textrm{otherwise}
\end{array} \right.
\end{displaymath}
Let $I' \subseteq I$ be a set of full Lebesgue measure for which \eqref{zlatos-decay-for-phi1} and 
\eqref{zlatos-decay-for-phi2} hold. Since the spectral measure for the unperturbed operator is 
$(1-2\beta_2)$-continuous for almost every boundary condition, it follows from the theory of rank-one perturbations, from 
the fact that $H_\theta(v,\gamma)$ has no absolutely continuous spectrum on $I'$ and from Theorem~1.4 in \cite{kls}, that 
there exists a set $I''\subseteq I'$ of full Lebesgue measure such that for any $E \in I''$
\beq \label{beta(E)-zlatos}
\beta(E) \geq \frac{1-2\beta_2}{2-(1-2\beta_2)} =\frac {1-2\beta_2}{1+2\beta_2}>0,
\eeq
and therefore
\beq \label{eta(E)-zlatos}
\eta(E) = \frac{1-\beta(E)}{\beta(E)}\leq \frac{4\beta_2}{1-2\beta_2}.
\eeq

 From the fact that almost every energy is regular, we may also assume that
\beq \label{Lambda0-intersect-J''}
\Lambda_0 \cap I'' = I''.
\eeq
Thus, in order to get almost sure stability of the asymptotic behavior of the generalized eigenfunctions on $I''$, we 
only need to show 
\beq \label{stability-condition-zlatos}
\sum_{n=1}^\infty r_{\ti{\eta}}^E(n) \mean{ \ti{b}_{s,\omega}(n)^2 }=
\sum_{n=1}^\infty r_{\ti{\eta}}^E(n) \frac{1}{n^{2s}}<\infty 
\eeq
for every $E \in I''$ and some $\ti{\eta}>\eta(E)$ (recall \eqref{definition-rnE}). 
Given \eqref{s-condition-zlatos}, \eqref{eta(E)-zlatos} and 
\eqref{zlatos-decay-for-phi1}-\eqref{zlatos-decay-for-phi2}, it is easy to verify that this is indeed the case.
Thus, it follows from Theorem~\ref{singular_stability} that up to a set of Lebesgue measure zero
\beq \no
I'' \cap\Sigma_{\alphac}\left(\se{b^{v,\gamma}}\right)=\Sigma_{\alphac}\left(\sse{b^{v,\gamma}}{\ti{b}^0_{s,\omega}}\right)
\eeq
\beq \no
I'' \cap\Sigma_{(1-\alpha){\rm s}}\left(\se{b^{v,\gamma}}\right)
=\Sigma_{(1-\alpha){\rm s}}\left(\sse{b^{v,\gamma}}{\ti{b}^0_{s,\omega}}\right)
\eeq
for a.e.\ $\omega$. 

Now, from the fact that the probability distribution of $\ti{b}_{s,\omega}(1)$ is absolutely continuous, it follows, 
using the theory of rank one perturbations \cite{rankone}, that for almost every realization of the
random perturbing potential, the spectral measure of $H_\omega(v,\gamma,s)$ is $\alpha$-continuous and 
$(1-\alpha)$-singular on J.
\end{proof}



\begin{thebibliography}{99}
\small
\bibitem{berez} J.~Berezanskii, {\it Expansions in Eigenfunctions of 
Selfadjoint Operators}, Transl.\ Math.\ Monographs, Vol.\ {\bf 17},
Amer.\ Math.\ Soc.,
Providence, RI, 1968.

\bibitem{christ-kiselev} M.~Christ and A.~Kiselev, {\it Absolutely continuous 
spectrum for one-dimensional Schr\"odinger operators with slowly decaying 
potentials: some optimal results}, J.\ Amer.\ Math.\ Soc.\ {\bf 11} (1998), 771--797.

\bibitem{christ-kiselev-b} M.~Christ and A.~Kiselev,
{\it WKB and spectral analysis of one-dimensional Schr\"odinger
operators with slowly varying potentials},
Commun.\ Math.\ Phys.\ {\bf 218} (2001), 245--262.

\bibitem{christ-kiselev-remling} M.~Christ, A.~Kiselev, and C.~Remling,
{\it The absolutely continuous spectrum of 
one-dimensional Schr\"odinger operators with decaying potentials},
Math.\ Res.\ Lett.\ {\bf 4} (1997), 719--723.

\bibitem{deift-killip} P.~Deift and R.~Killip, {\it On the absolutely continuous 
spectrum of one dimensional Schr\"odinger operators with square summable 
potentials}, Commun.\ Math.\ Phys.\ {\bf 203} (1999), 341--347.

\bibitem{delrio} R.~Del-Rio, B.~Simon, and G.~Stolz, {\it Stability of spectral 
types for Sturm-Liouville operators}, Math.\ Res.\ Lett.\ {\bf 1} (1994), 437--450.

\bibitem{djls} R.~Del-Rio, S.~Jitomirskaya, Y.~Last, and B.~Simon, {\it Operators 
with singular continuous spectrum, {\rm IV}.\ Hausdorff dimensions, rank one 
perturbations, and localization}, J.\ d'Analyse Math.\ {\bf 69} (1996), 153--200.

\bibitem{subord} D.J.~Gilbert and D.B.~Pearson,
{\it On subordinacy and analysis of the 
spectrum of one-dimensional Schr\"odinger operators},
J.\ Math.\ Anal.\ Appl.\ {\bf 128} (1987), 30--56.

\bibitem{jit-last} S.~Jitomirskaya and Y.~Last, {\it Power-law subordinacy and 
singular spectra, {\rm I}.\ Half-line operators},
Acta Math.\ {\bf 183} (1999), 171--189.

\bibitem{kala} U.~Kaluzhny and Y.~Last, {\it Purely absolutely continuous 
spectrum for some random Jacobi matrices},
Proceedings of ``Probability and Mathematical Physics'' a conference 
in honor of Stanislav Molchanov's 65'th birthday, to appear.

\bibitem{killip} R.~Killip, {\it Perturbations of one-dimensional Schr\"odinger 
operators preserving the absolutely continuous spectrum},
Int.\ Math.\ Res.\ Not.\ {\bf 38} (2002), 2029--2061.

\bibitem{killip-simon} R.~Killip and B.~Simon, {\it Sum rules for Jacobi matrices 
and their applications to spectral theory},
Ann.\ Math.\ {\bf 158} (2003), 253--321.

\bibitem{kiselev1} A.~Kiselev, {\it Absolutely continuous spectrum of 
one-dimensional Schr\"odinger operators and Jacobi matrices with slowly 
decreasing potentials}, Commun.\ Math.\ Phys.\ {\bf 179} (1996), 377--400.

\bibitem{kiselev2} A.~Kiselev, {\it Stability of the absolutely continuous 
spectrum of the Schr\"odinger equation under slowly decaying perturbations 
and a.e.\ convergence of integral operators},
Duke Math.\ J.\ {\bf 94} (1998), 619--646.

\bibitem{efgp} A.~Kiselev, Y.~Last, and B.~Simon, {\it Modified Pr\"ufer and EFGP 
transforms and the spectral analysis of one-dimensional Schr\"odinger 
operators}, Commun.\ Math.\ Phys.\ {\bf 194} (1998), 1--45. 

\bibitem{kls} A.~Kiselev, Y.~Last, and B.~Simon, {\it Stability of 
Singular Spectral types under decaying perturbations}, 
J.\ Funct.\ Anal.\ {\bf 198} (2003), 1--27.

\bibitem{hausdorff} Y.~Last, {\it Quantum dynamics and decompositions of 
singular continuous spectra}, J.\ Funct.\ Anal.\ {\bf 142} (1996), 406--445.

\bibitem{last-simon} Y.~Last and B.~Simon, {\it Eigenfunctions, transfer matrices, 
and absolutely continuous spectrum of one-dimensional Schr\"odinger 
operators}, Invent.\ Math.\ {\bf 135} (1999), 329--367.

\bibitem{schrodinger-conjecture} V.~P.~Maslov, S.~A.~Molchanov, and A.~Ya.~Gordon, 
{\it Behavior of generalized eigenfunctions at infinity and the 
Schr\"odinger conjecture}, Russian J.\ Math.\ Phys.\ {\bf 1} (1993), 71--104.

\bibitem{reed-simon1} M.~Reed and B.~Simon, 
{\it Methods of Modern Mathematical Physics, {\rm I}.\ Functional Analysis}, 
Academic Press, New York, 1972.

\bibitem{reed-simon3} M.~Reed and B.~Simon,
{\it Methods of Modern Mathematical Physics, {\rm III}.\ Scattering Theory},
Academic Press, New York, 1979.

\bibitem{remling} C.~Remling, {\it The absolutely continuous spectrum of 
one-dimensional Schr\"odinger operators with decaying potentials},
Commun.\ Math.\ Phys.\ {\bf 193} (1998), 151--170.

\bibitem{rankone} B.~Simon, {\it Spectral analysis of rank one perturbations 
and applications}, in ``Proc.\ Mathematical Quantum Theory, II.\ Schr\"odinger Operators''
(Vancouver, Canada, 1993), 
pp.\ 109--149, CRM Proceedings and Lecture Notes, {\bf 8},
American Mathematical Society, Providence, RI, 1995.

\bibitem{bounded-eig} B.~Simon, {\it Bounded eigenfunctions and absolutely 
continuous spectra for one-dimensional Schr\"odinger operators}, Proc.\ 
Amer.\ Math.\ Soc.\ {\bf 124} (1996), 3361--3369.

\bibitem{zlatos} A.~Zlato\v s, {\it Sparse potentials with fractional Hausdorff 
dimensions}, J.\ Funct.\ Anal.\ {\bf 207} (2004), 216--252.

\end{thebibliography}
\end{document}